\documentclass[MSNbibl,number,citesort,dvips]{arxbj}
\usepackage{stfloats}
\usepackage{graphicx}

\aid{0}
\volume{19}
\issue{2}
\pubyear{2013}
\firstpage{426}
\lastpage{461}
\doi{10.3150/11-BEJ407}

\makeatletter
\fnbelowfloat
\newcommand{\eqref}[1]{(\ref{#1})}
\newcommand{\bea}{\begin{equation}}
\newcommand{\eea}{\end{equation}}
\newcommand{\calf}{\mathcal{F}}
\newcommand{\calg}{\mathcal{G}}
\newcommand{\ep}{\varepsilon}
\newcommand{\R}{\mathbb{R}}
\newcommand{\bbR}{\mathbb{R}}
\newcommand{\bbU}{\mathbb{U}}

\newcommand{\E}{\mathbb{E}}
\newcommand{\PP}{\mathbb{P}}
\newcommand{\F}{\mathbb{F}}
\newcommand{\circj}{\widetilde {J}^n}
\newcommand{\circy}{\widetilde{Y}}

\newtheorem{thmm}{Theorem}
\newproclaim{definition}{Definition}
\newtheorem{lemm}{Lemma}
\newproclaim{ass}{Assumption}
\newtheorem{propm}{Proposition}
\newremark{rem}{Remark}
\newremark{notation}{Notation}
\makeatother

\begin{document}
\begin{frontmatter}

\title{Estimation of the lead-lag parameter from non-synchronous data}
\runtitle{Estimation of the lead-lag parameter}

\begin{aug}
%%%% inicialai - be tarpu
\author[a]{\fnms{M.} \snm{Hoffmann}\thanksref{a}\ead[label=e1]{marc.hoffmann@ensae.fr}},
\author[b]{\fnms{M.} \snm{Rosenbaum}\corref{}\thanksref{b}\ead[label=e2]{mathieu.rosenbaum@polytechnique.edu}}
\and
\author[c]{\fnms{N.} \snm{Yoshida}\thanksref{c}\ead[label=e3]{nakahiro@ms.u-tokyo.ac.jp}}
\runauthor{M. Hoffmann, M. Rosenbaum and N. Yoshida} %% auto
\address[a]{ENSAE -- CREST and CNRS UMR 8050, Timbre J120, 3, avenue Pierre Larousse,
92245 Malakoff Cedex, France. \printead{e1}}
\address[b]{LPMA -- Universit\'{e} Pierre et Marie Curie (Paris 6) and CREST, 4 Place Jussieu,
75252 Paris Cedex 05, France.
\printead{e2}}
\address[c]{University of Tokyo and Japan Science and Technology
Agency, Graduate School of Mathematical Sciences, University of Tokyo, 3-8-1 Komaba, Meguro-ku,
Tokyo 153-8914, Japan.\\ \printead{e3}}
\end{aug}

% HISTORY:
\received{\smonth{1} \syear{2010}}
\revised{\smonth{9} \syear{2011}}

% ABSTRACT
%
\begin{abstract}
We propose a simple continuous time model for modeling the lead-lag
effect between two financial assets. A two-dimensional process
$(X_t,Y_t)$ reproduces a lead-lag effect if, for some time shift
$\vartheta\in\R$, the process $(X_t,Y_{t+\vartheta})$ is a
semi-martingale with respect to a certain filtration. The value of
the time shift $\vartheta$ is the lead-lag parameter. Depending on
the underlying filtration, the standard no-arbitrage case is
obtained for $\vartheta= 0$. We study the problem of estimating the
unknown parameter $\vartheta\in\R$, given randomly sampled
non-synchronous data from $(X_t)$ and $(Y_t)$. By applying a certain
contrast optimization based on a modified version of the
Hayashi--Yoshida covariation estimator, we obtain a consistent
estimator of the lead-lag parameter, together with an explicit rate
of convergence governed by the sparsity of the sampling design.
\end{abstract}

% KEYWORDS
%

\end{frontmatter}

%s1 ###
\section{Introduction}
Market participants usually agree that
certain pairs of assets $(X,Y)$ share a ``lead-lag effect,'' in the
sense that the lagger (or follower)
price process $Y$ tends to partially reproduce the oscillations of
the leader (or driver) price process $X$, with some temporal
delay, or vice-versa. This property is usually referred to as the
``lead-lag effect.'' The lead-lag effect may have some importance in practice,
when assessing the quality of risk management indicators, for
instance, or, more generally, when considering statistical arbitrage
strategies. Also, note that it can be measured at various temporal
scales (daily,
hourly or even at the level of seconds, for flow products traded
on electronic markets).

The lead-lag effect is a concept of common practice that has some
history in
financial econometrics. In time series for instance, this notion
can be linked to the concept of Granger causality, and we refer to
Comte and Renault~\cite{CR} for a general approach. From a
phenomenological perspective, the lead-lag effect is supported by empirical
evidence reported in~\cite{empi1,empi2} and~\cite{empi3},
together with~\cite{O} and the references therein. To our
knowledge, however, only few mathematical results are available
from the point of view of statistical estimation from discretely
observed, continuous-time processes. The purpose of this paper is
to -- partly -- fill in this gap. (Also, recently, Robert and Rosenbaum
study in~\cite{RR3} the lead-lag effect by means of random matrices,
in a
mixed asymptotic framework, a setting which is relatively different
than in the present paper.)
% We will work here in a usual high
%frequency setting.
%ee together with the references therein.

%s1.1 ###
\subsection{Motivation}
\begin{longlist}[(1)]
%In this paper, we investigate the lead-lag effect from two
%complementary angles:

\item[(1)] Our primary goal is to provide a simple -- yet
relatively general -- model for capturing the lead-lag effect in
continuous time, readily compatible with stochastic calculus in
financial modeling.
%, \textit{i.e.} the presence of a lead-lag temporal parameter $\vartheta
Informally, if $\tau_{-\vartheta}(Y)_t:=Y_{t+\vartheta}$, with
$\vartheta\in\R$, is the time-shift operator, we say that the pair
$(X,Y)$ will produce a lead-lag effect as soon as $(X,\tau
_{-\vartheta}(Y))$ is a (regular) semi-martingale with respect to
an appropriate filtration, for some $\vartheta$, called the lead-lag
parameter. The usual no-arbitrage case is embedded into this framework
for $\vartheta=0$. More in Section~\ref{leadlag model} below.

\item[(2)] At a similar level of importance, we aim at
constructing a simple and efficient procedure for estimating the
lead-lag parameter $\vartheta$ based on historical data. The
underlying statistical model is generated by a -- possibly random
-- sampling of both $X$ and $Y$. The sampling typically happens at
irregularly and non-synchronous times for $X$ and $Y$. We
construct, in the paper, an estimator of $\vartheta$ based on a
modification of the Hayashi--Yoshida covariation estimator; see
\cite{HY1} and~\cite{HY2}. Our result is that the lead-lag
parameter can be consistently estimated against a fairly general
class of sampling
schemes. Moreover, we explicit the rate of convergence of our
procedure.

%If $\{s_i\}$ and $\{t_i\}$ denote the sampling times of $X$ and $Y$
%respectively, we find a rate of convergence which is essentially
%governed by $\max\{\max(s_i-s_{i-1}), \max\{t_i-t_{i-1}\}\}$. In
%particular, if we have regular sampling schemes over the time horizon

%Addressing the modelling and the statistical identification of the LLE
%parameter raises two kinds of intertwined difficulties:

\item[(3)] From a financial point of view, unless
appropriate time shifts are operated, our model incapacitates our
primary assets $X$ and $Y$ to be a semi-martingale with respect to the
same filtration. This is consistent, as far as modeling is concerned,
but allows, in principle, for market imperfections such as statistical
arbitrage if the lead-lag parameter $\vartheta$ is different from zero.
More in Section~\ref{discussion} below. Addressing such a possibility
is indeed the issue of the lead-lag effect, but we will content
ourselves with detecting whether the lead-lag effect is present or not.
The quantization of statistical arbitrage in terms of $\vartheta$ (and
other parameters such as trading frequency, market friction, volatility
and so on) lies beyond the scope of this paper.
%Hopefully, this will be set out in details elsewhere.
%%\item From a financial modelling point of view,
%arbitrage

\item[(4)] From a statistical inference point of view, the
statistician and the data provider are not necessarily the same
agents, and this leads to technical difficulties linked to the
sampling strategy. The data provider may choose the
opening/closing for $X$ and $Y$, possibly traded on different
markets, possibly on different time clocks. He or she may also
sample points at certain trading times or events which are
randomly chosen in a particular time window. This typically
happens if daily data are considered. At a completely different
level, if high-frequency data are concerned, trading times are
genuinely random and non-synchronous. Our approach will
simultaneously incorporate these different points of view.
\end{longlist}
%s1.2 ###
\subsection{Organization of the paper}
In Section~\ref{leadlag model}, we present our stochastic model for
describing the lead-lag effect. We start with the simplest Bachelier
model with no drift in Section~\ref{leadlag Bachelier}. The issue
boils down to defining properly the lead-lag effect between two
correlated Brownian motions. In Section~\ref{leadlag general}, a
general lead-lag model is presented for two-dimensional process, for
which the marginal processes are semi-martingales with locally
bounded drift and continuous local martingale part, with properly
defined diffusion coefficients.

We present our main result in Section~\ref{main result}. Section
\ref{statmodel} gives a precise construction of the underlying
statistical experiment with the corresponding assumptions on the
observation sampling schemes. The estimation procedure is
constructed in Section~\ref{procedure}, via an appropriate contrast
function based on the covariation between $X$ and $Y$ when one asset is
artificially shifted in time, the amount of this shift being the
argument of the contrast function. Our estimator is robust to
non-synchronous data and does not require any pre-processing contrary
to the previous tick algorithm; see, for example, \cite
{zhang2010estimating}. In Section
\ref{th}, we state our main result in Theorem~\ref{mainth}: we show
that the lead-lag parameter
between $X$ and $Y$ can be consistently estimated from non-synchronous
historical data over a fixed time horizon
$[0,T]$. The rate is governed by $\Delta_n$, the
maximal distance between two data points.
%In particular, if the
%general setting is reduced to regular and synchronous sampling
%frequency $\Delta^{-1}$ between two data points, the rate of
We show that the rate of convergence of our estimator is essentially
$\Delta_n^{-1}$ and not
$\Delta_n^{-1/2}$, as one would expect from a regular estimation
problem in diffusion processes; see, for example,~\cite{GCJ}.
This comes from the underlying structure of the statistical model,
which is not regular, and which shares some analogy with change-point
problems. As for our procedure, we investigate further its asymptotic
properties in Proposition~\ref{asympt_behaviour_contrast} when we
confine ourselves to the simpler case where $X$ and $Y$ are marginally
Brownian motions that are observed at synchronous data points. In that
case, we can exhibit a central limit theorem for our contrast function.
A closer inspection of the limiting variance reveals the effect of the
correlation between the two assets, which also plays a role in the
accuracy of the estimation procedure. Finally, we show in Proposition
\ref{noclt} that a simple central limit theorem cannot hold for our
estimator. We discuss this effect which is somewhat linked to the
discretisation of our method.

Theorem~\ref{mainth} is good news, as far as practical implementation
is concerned, and is further addressed in the discussion in Section
\ref
{discussion}, appended with numerical illustrations on simulated data
in Section~\ref{illustration} and on real data in Section \ref
{illustration real data}. The proofs are delayed until Section \ref
{proofs} and the \hyperref[appendix]{Appendix} contains auxiliary
technical results.

%s2 ###
\section{The lead-lag model} \label{leadlag model}
%s2.1 ###
\subsection{The Bachelier model} \label{leadlag Bachelier}
A simple lead-lag Bachelier model with no drift between two Brownian
motion components can be described as follows. On a filtered space
$(\Omega, {\mathcal F}, \mathbb{F} = ({\mathcal F}_t)_{t \geq0},\PP
)$, we consider a two-dimensional $\mathbb{F}$-Brownian motion
$B=(B^{(1)}, B^{(2)})$ such that $\langle B^{(1)},B^{(2)}\rangle
_t=\rho
  t$ for every $t \geq0$ and for some $\rho\in[-1,1]$. Let $T>0$\vadjust{\goodbreak} be
some terminal time, fixed throughout the paper.
For $t \in[0,T]$, set
%We start with a random process $(X,Y)$ with value in $\R^2$ by
%setting, for $t\in[0,T]$,
%
\[
\cases{
X_t  :=  x_0+\sigma_1 B^{(1)}_t, \vspace*{2pt}\cr
\widetilde Y_t  :=  y_0+\sigma_2 B^{(2)}_t,}
\]
%for a two-dimensional Bachelier model with no drift, with
where $x_0,y_0 \in\R$ and $\sigma_1>0$, $\sigma_2 >0$ are
given constants.
%This is a Bachelier type model with no drift. The process
%and is a two-dimensional Brownian motion
%for some filtration $\mathbb{F}=(\mathcal{F}_t)_{t\geq0}$.
%We
%further assume that
%for some $\rho\in[0,1]$.
The corresponding Black--Scholes version
of this model is readily obtained by exponentiating $X$ and $\widetilde
Y$. We introduce a lead-lag effect between $X$ and
$\widetilde Y$ by operating a time shift: let $\vartheta\in\R$ represent
the lead or lag time between $X$ and $\widetilde Y$ (and assume for
simplicity that
$\vartheta\geq0$). Put
\begin{equation} \label{shiftoperator}
\tau_\vartheta(\widetilde Y)_t:=\widetilde Y_{t -\vartheta},\qquad  t
\in
[\vartheta,T].
\end{equation}
%
%where $\vartheta^+=\max\{\vartheta,0\}$ and $\vartheta^-=-\min\{
Our lead-lag model is the two-dimensional process
\[
(X,\tau_\vartheta(\widetilde Y))=(X_t,\tau_\vartheta
(\widetilde Y)_t)_{t \in[\vartheta,T]}.
\]
Since we have
$B^{(2)}_{t} =\rho B^{(1)}_{t}+(1-\rho^2)^{1/2} W_{t}$
with $W=(W_t)_{t \in[0,T]}$, a Brownian motion independent of $B^{(1)}$,
we obtain the simple and explicit representation
\begin{equation} \label{leadlagrepresentation}
\cases{
X_t  =  x_0+\sigma_1 B^{(1)}_t,\vspace*{2pt}\cr
\tau_\vartheta(\widetilde Y)_{t}  =  y_0+ \rho \sigma_2
B^{(1)}_{t-\vartheta}+\sigma_2(1-\rho^2)^{1/2} W_{t-\vartheta}}
\end{equation}
for $t \in[\vartheta,T]$. In this representation, the interpretation
of the lead-lag parameter $\vartheta$ is transparent. Alternatively, if
we start with a process $(X,Y)$ having representation
\begin{equation} \label{defquickbachelier}
(X,Y) = (X,\tau_\vartheta(\widetilde Y))
\end{equation}
as in \eqref{leadlagrepresentation}, the lead-lag interpretation
between $X$ and $Y$ readily follows. Since
$\vartheta
\geq0$, the sample path of $X$ anticipates on the path of $Y$ by a
time shift $\vartheta$ and to an amount -- measured in
normalized standard deviation -- proportional to $\rho  \sigma
_2/\sigma_1$.
In that case, we say that $ X$ is the leader, and $Y$ is
the lagger. For the case $\vartheta< 0$, we intertwine the roles of $ X$
and $Y$ in the terminology.
\begin{rem}
Note that, except in the case $\vartheta=0$, the process
$(X_t,Y_t)_{t \in[\vartheta,T]}$ is not an
$\mathbb{F}$-martingale. However, each component is a martingale
with respect to a different filtration: $X$ is an
$\mathbb{F}$-martingale and $Y = \tau_{\vartheta}(\widetilde Y)$ is an
$\mathbb{F}^{\vartheta}$-martingale, with
$\mathbb{F}^{\vartheta}=(\mathcal{F}^\vartheta_t)_{t\geq
\vartheta}$, and
\mbox{$\mathcal{F}^{\vartheta}_t=\mathcal{F}_{t-\vartheta}$}.
%In the representation \eqref{defquickbachelier}, the process $\tau_
%the ``point of view'' of $\mathbb{F}$.
%: the filtration generated by $\tau_\vartheta(\widetilde Y)$ differs
%from $\mathbb{F}$.
\end{rem}
%
%s2.2 ###
\subsection{Lead-lag between two semi-martingales} \label{leadlag general}
We generalize the lead-lag model \eqref{defquickbachelier} to
semi-martingales with local martingale components that can be
represented as It\^o local martingales.

We need some notation. Let $T>0$ be some terminal time, and let $\delta
>0$ represent the maximum temporal lead-lag allowed for the model,
fixed throughout the paper. On a probability space $(\Omega
,{\mathcal F}, \PP)$, let
$\mathbb{F}= ({\mathcal F}_t)_{t \in[-\delta, T+\delta]}$
be a
filtration
satisfying the usual conditions. We denote by
$\mathbb{F}_{[a,b]} = ({\mathcal F}_t)_{t \in[a,b]}$
the restriction of $\mathbb{F}$ to the time interval $[a,b]$.\vadjust{\goodbreak}

\begin{definition}\label{regmart}
The two-dimensional process $(X,Y)_{t \in[0,T+\delta]}$ is a
regular semi-martingale with lead-lag parameter $\vartheta\in
[0,\delta)$ if the following decomposition holds:
\[
X=X^c+A,\qquad  Y=Y^c+B,
\]
with the following properties:
\begin{itemize}
\item The process $(X^c_t)_{t\in[0,T+\delta]}$ is a continuous
$\mathbb
{F}_{[0,T+\delta]}$-local martingale,
and the process $(Y^c_t)_{t\in[0,T+\delta]}$ is a continuous $\mathbb
{F}^{\vartheta}_{[0,T+\delta]}$-local martingale.
\item The quadratic variations $\langle X^c \rangle_{t\in[0,T+\delta]}$
and $\langle Y^c \rangle_{t\in[0,T+\delta]}$ are absolutely continuous
w.r.t. the Lebesgue measure, and their Radon--Nikodym derivatives
admit a locally bounded version.
\item The drifts $A$ and $B$ have finite variation over
$[0,T+\delta]$.\vspace*{-1pt}
\end{itemize}
\end{definition}

\begin{definition}
The two-dimensional process $(X,Y)_{t \in[0,T+\delta]}$ is a
regular semi-martingale with lead-lag parameter $\vartheta\in
(-\delta,0]$ if the same properties as in Definition~\ref{regmart}
hold, with $X$ and $Y$ intertwined and $\vartheta$ replaced by
$-\vartheta$.\vspace*{-1pt}
\end{definition}

\begin{rem}
If $(X,Y)_{t \in[0,T+\delta]}$ is a regular semi-martingale with
lead-lag parameter $\vartheta\in[0,\delta)$, then the process
$(\tau_{-\vartheta}(Y^c))_{t\in[-\vartheta,T]}$ is a continuous
$\mathbb{F}_{[-\vartheta,T]}$-local martingale, with
$\tau_{-\vartheta}(Y)_t=Y_{t+\vartheta}$ the (inverse of the) shift
operator defined in \eqref{shiftoperator}.\vspace*{-1pt}
\end{rem}

\begin{rem}\label{inverse}
If $(X,Y)_{t \in[0,T+\delta]}$ is a regular semi-martingale with
lead-lag parameter $\vartheta\in[0,\delta)$, then the process
$(Y,X)_{t \in[0,T+\delta]}$ is a regular semi-martingale with
lead-lag parameter $-\vartheta$.\vspace*{-1pt}
\end{rem}

%
%, giving precise definition in Section~\ref{defgenleadlag}. We revert
%the scheme presented in the Bachelier case. Let
%$$(X,Y) = (X_t,Y_t)_{t \in{\mathcal I}},  {\mathcal I}\subset\R$$
%be a two-dimensional process,
%where ${\mathcal I}$ is an interval, and such that both $X$ and $Y$
%can reasonably be extended over ${\mathcal I}_\vartheta= \{t, t\pm
%lead-lag representation property with lead-lag parameter $\vartheta$
%if $(X,\tau_\vartheta(Y))$ is a regular semimartingale.
%model is a two-dimensional price process that can be mapped into a
%regular semimartingale by an appropriate time shift on one of its
%components.

%s3 ###
\section{Main result} \label{main result}\vspace*{-1pt}
%s3.1 ###
\subsection{The statistical model} \label{statmodel}\vspace*{-1pt}

We observe a two-dimensional price process $(X,Y)$
at discrete times. The components $X$ and $Y$ are observed over
the time horizon $[0,T+\delta]$. The following assumption is in force
throughout:\vspace*{-1pt}

\renewcommand{\theass}{\Alph{ass}}
\begin{ass}\label{assa}
The process $(X,Y) = (X_t,Y_t)_{t \in[0,T+\delta]}$ is a regular
semi-martingale with lead-lag parameter $\Theta\in\vartheta=
(-\delta
, \delta)$.
\end{ass}

The -- possibly random -- observation times are given by the following
subdivisions of $[0,T+\delta]$:
\begin{equation} \label{subdivisionX}
{\mathcal T}^X:=\{s_{1,n_1} < s_{2,n_1} < \cdots<
s_{n_1,n_1}\}
\end{equation}
for $X$ and
\begin{equation} \label{subdivisionY}
{\mathcal T}^Y:=\{t_{1,n_2} < t_{2,n_2} < \cdots<
t_{n_2,n_2}\}
\end{equation}
for $Y$, with $n_1=n_2$ or not. For simplicity, we assume
$s_{1,n_1}=t_{1,n_2}=0$ and $s_{n_1,n_1}=t_{n_2,n_2}=T+\delta$.
The sample points are either chosen by the statistician or\vadjust{\goodbreak}
dictated for practical convenience by the data provider. They are
usually neither equispaced in time nor synchronous, and may depend
on the values of $X$ and $Y$.
% in a
%certain manner, see the Introduction.

For some unknown $\vartheta\in\Theta:=(-\delta,\delta)$, the
process $(X,Y)$ is a regular semi-martingale with lead-lag parameter
$\vartheta$, and we want to estimate $\vartheta$ based on the set of
historical data
%$$\{X_{s_{i,n_1}}, i=1,\ldots, n_1\} \cup
%$$
%
\begin{equation} \label{data set}
\{X_{s}, s\in{\mathcal T}^X\} \cup\{Y_{t}, t\in
{\mathcal T}^Y\}.
\end{equation}
%
%that generate a statistical experiment with parameter $\vartheta$
%and parameter space
%$\vartheta\in\vartheta=(-\delta, \delta)$.
In order to describe precisely the property of the sampling scheme
${\mathcal T}^X \cup{\mathcal T}^Y$, we need some notation that we
borrow from Hayashi and Yoshida~\cite{HY1}. The subdivision ${\mathcal T}^X$
% and ${\mathcal T}^Y$
introduced in \eqref{subdivisionX}
% and \eqref{subdivisionY}
is mapped into a family of intervals
\begin{equation} \label{defI}
{\mathcal I} = \{I=(\underline{I}, \overline{I}] = (s_{i,n_1},
s_{i+1,n_1}],  i=1,\ldots, n_1-1\}.
\end{equation}
Likewise, the subdivision ${\mathcal T}^Y$ defined in \eqref
{subdivisionY} is mapped into
\[
{\mathcal J} = \{J=(\underline{J}, \overline{J}] =
(t_{j,n_2}, s_{j+1,n_2}],  j=1,\ldots, n_2-1\}.
\]
We will
systematically employ the notation $I$ (resp., $J$) for an element
of ${\mathcal I}$ (resp., ${\mathcal J}$). We set
\[
\Delta_n:=\max\bigl\{\sup\{|I|, I \in{\mathcal I}\}, \sup\{
|J|, J
\in{\mathcal J}\}\bigr\},
\]
where $|I|$ (resp., $|J|$) denotes the length of the interval $I$
(resp.,
$J$), and $n$ is a parameter tending to infinity. %Indeed, we will see
%that the (random) quantity $\Delta_n$ is relevant for measuring the
%accuracy of estimation of the lead-lag parameter.

\begin{rem}
One may think of $n$ being the number of data points extracted from the
sampling, that is, $n=\sharp{\mathcal I}+\sharp{\mathcal J}$.
However, as we will see, only the (random) quantity $\Delta_n$ will
prove relevant for measuring the accuracy of estimation of the lead-lag
parameter.
%, so we may forget about the original interpretation of $n$ and let it
%be an arbitrary asymptotic parameter. This will allow for instance to
%include in our setting more general sampling schemes.
\end{rem}

The assumptions on the sampling scheme is the following.

\begin{ass}\label{assb}
\begin{enumerate}[B1.]
\item[B1.] There exists a deterministic sequence of positive
numbers $v_n$
such that $v_n<\delta$ and $v_n\rightarrow0$ as $n\rightarrow
\infty$. Moreover
\[
v_n^{-1} \Delta_n\rightarrow0
\]
in probability as $n \rightarrow\infty$.
%$v_n/r_n\stackrel{\PP}{\longrightarrow}0,$
%where $\stackrel{\PP}{\rightarrow}$ denotes convergence in probability.
%
\item[B2.] For all $I \in{\mathcal I}$, the random times
$\underline{I}$ and $\overline{I}$ are $\mathbb{F}^{v_n}$-stopping
times if $\vartheta\geq0$ (resp., $\mathbb{F}^{-\vartheta
+v_n}$-stopping times if $\vartheta<0$). For all $J \in{\mathcal J}$,
the random times $\underline{J}$ and $\overline{J}$ are $\mathbb
{F}^{\vartheta+v_n}$-stopping times if $\vartheta\geq0$ (resp.,
$\mathbb{F}^{v_n}$-stopping times if $\vartheta<0$).
% where $\mathbb{F}^{\vartheta+v_n}$ is defined through
%
\item[B3.] There exists a finite grid ${\mathcal G}^n \subset
\Theta$ such that $0 \in{\mathcal G}^n$ and
\begin{itemize}[--]
\item[--] For some $\gamma>0$, we have $\sharp  {\mathcal G}^n =
\mathrm{O}(v_n^{-\gamma})$.
\item[--] For some deterministic sequence $\rho_n >0$, we have
\[
\bigcup_{\widetilde\vartheta\in{\mathcal G}^n}[\widetilde\vartheta
-\rho_n,\widetilde\vartheta+\rho_n] \supset\Theta
\]
and
\[
\lim_{n \rightarrow\infty}\rho_n \min\{\E[\sharp
{\mathcal
I}], \E[\sharp  {\mathcal J}]\} \rightarrow0.
\]
%
%as $n\rightarrow\infty$.
\end{itemize}
\end{enumerate}
\end{ass}

\begin{rem}
Since both $\E[\sharp{\mathcal I}]$ and $\E[\sharp
{\mathcal J}]$ diverge at rate no less than $v_n^{-1}$, Assumption~B3
implies that $\rho_n = \mathrm{o}(v_n)$. With no loss of generality, we thus
may (and will) assume that $\rho_n \leq v_n$ for all $n$.
\end{rem}
%
%As will appear in Theorem~\ref{} %below, the order of magnitude of the
%number of data $n_1+n_2$ is irrelevant, provided ${\mathcal I}$ and ${
%assumption A.
%{\tt[note on strong predictability here]}
%The following assumptions on the lead-lag model $(X,Y)$ are in force
%throughout:\\

% \textbf{Assumption A.}\\
% \textbf{Assumption B.}\\

%s3.2 ###
\subsection{The estimation procedure}\label{procedure}
\subsubsection*{Preliminaries}
Assume first that the data arrive at regular and synchronous time
stamps over the time interval $[0,T] = [0,1]$, with $\Delta_n=1/n$ for
simplicity.
This means that we have $2n+2$ observations
\[
(X_0,Y_0), (X_{1/n}, Y_{1/n}), (X_{2/n}, Y_{2/n}),\ldots, (X_{1}, Y_{1}).
\]
For every integer $k \in\mathbb{Z}$, we form the shifted time series
\[
Y_{(k+i)/n}, \qquad i=1,2,\ldots
\]
for every $i$ such that $(k+i)/n$ is an admissible time stamp\footnote
{Possibly, we end up with an empty data set.}.
%For a continuous time process $Z$, we set $\delta_{i,\Delta_n}Z=:Z_{i
We can then construct the empirical covariation estimator
%{\mathcal C}_n(k):= & \sum_{i} (\delta_{i, \Delta_n}X)  (
%= &\sum_i (X_{i\Delta_n}-X_{(i-1)\Delta_n})(Y_{(i+k)
%
\[
{\mathcal C}_n(k):= \sum_i \bigl(X_{i/n}-X_{(i-1)/n}\bigr)
\bigl(Y_{(i+k)/n}-Y_{(i+k-1)/n}\bigr),
\]
where the sum in $i$ expands over all relevant data points. Over the
time interval $[0,1]$, the number of elements used for the computation
of ${\mathcal C}_n(k)$ should be of order $n$ as $n \rightarrow\infty
$. Assume further for simplicity that the process
$(X,Y)$ is a lead-lag Bachelier model in the sense of Section \ref
{leadlag Bachelier}, with lead-lag parameter $\vartheta=\vartheta
_n=k^0_n/n$, with $k^0_n$ an integer. On the one hand, for $k=k^0_n$,
we have the decomposition
\[
{\mathcal C}_n(k^0_n)=T_n^{(1)}+T_n^{(2)},
\]
with
\begin{eqnarray*}
T_n^{(1)}&=&\rho \sigma_1\sigma_2\sum_i
\bigl(B^{(1)}_{i/n}-B^{(1)}_{(i-1)/n}\bigr)^2,\\
T_n^{(2)}&=&\sqrt{1-\rho^2}\sigma_1\sigma_2\sum_i
\bigl(B^{(1)}_{i/n}-B^{(1)}_{(i-1)/n}\bigr)\bigl(W_{i/n}-W_{(i-1)/n}\bigr).
\end{eqnarray*}
Computing successively the fourth-order moment of the random variables
$T_n^{(1)}-\rho \sigma_1\sigma_2$ and $T_n^{(2)}$ and applying
Markov's inequality and the Borel--Cantelli lemma,\vadjust{\goodbreak} elementary
computations show
that $T_n^{(1)}\rightarrow\rho \sigma_1\sigma_2$ and
$T_n^{(2)}\rightarrow0$ as $n \rightarrow\infty$ almost surely, and
we derive
\[
{\mathcal C}_n(k^0_n)\rightarrow\rho \sigma_1\sigma_2 \qquad  \mbox
{as }
  n \rightarrow\infty  \mbox{ almost surely}.
\]
On the other hand, for $k\neq k^0_n$, we have
\[
{\mathcal C}_n(k)=\widetilde{T}_n^{(1)}+\widetilde{T}_n^{(2)},
\]
with
\begin{eqnarray*}
\widetilde{T}_n^{(1)}&=&\rho \sigma_1\sigma_2\sum_i
\bigl(B^{(1)}_{i/n}-B^{(1)}_{(i-1)/n}\bigr)
\bigl(B^{(1)}_{(i+k-k^0_n)/n}-B^{(1)}_{(i+k-k^0_n-1)/n}\bigr),\\
\widetilde{T}_n^{(2)}&=&\sqrt{1-\rho^2}\sigma_1\sigma_2\sum_i
\bigl(B^{(1)}_{i/n}-B^{(1)}_{(i-1)/n}\bigr)
\bigl(W_{(i+k-k^0_n)/n}-W_{(i+k-k^0_n-1)/n}\bigr).
\end{eqnarray*}
Thus, for fixed $n$ and $k>k^0_n$, the process
\[
j \leadsto\sum_{i=1}^j \bigl(X_{i/n}-X_{(i-1)/n}\bigr)
\bigl(Y_{(i+k)/n}-Y_{(i+k-1)/n}\bigr)
\]
is $(\mathcal{F}_{(j+k-k^0_n)/n})_{j\geq1}$-martingale. Consequently,
using the
Burkholder--Davis--Gundy inequality, we easily obtain that
\[
\E[{\mathcal C}_n(k)^6]\leq c n^{-3},
\]
up to some constant $c >0$. The same result holds for $k<k^0_n$. We infer
\[
\E\Bigl[\Bigl(\mathop{\operatorname{sup}}_{k\neq k^0_n}|{\mathcal
C}_n(k)|
\Bigr)^6\Bigr]\leq c n^{-2}
\]
up to a modification of $c$. Using again Markov's inequality and the
Borel--Cantelli lemma, we finally obtain that
\[
\mathop{\operatorname{sup}}_{k\neq k^0_n}|{\mathcal C}_n(k)|\rightarrow0
\qquad
\mbox{as }   n\rightarrow\infty   \mbox{ almost surely}.
\]
%
%Heuristically, we have
%$${\mathcal C}_n(k) \approx\Delta_n^{-1}\E[(X_\cdot-X_{\cdot-
%$$
%{\mathcal C}_n(k) \approx
%as $\Delta_n \rightarrow0$, where the noise remainder
%term $\xi^n$ is a tight sequence by the central limit theorem.
%Now, thanks to the representation \eqref{leadlagrepresentation} we
%readily derive
%$$
% \Delta_n^{-1}\E[(\delta_{\cdot,\Delta_n} X)(\delta_{
%=
%0 & \operatorname{if} & k\neq k_0\\
%.$$
Therefore, provided $\rho\sigma_1\sigma_2\neq0$, we can detect
asymptotically the value $k^0_n$
that defines $\vartheta$ in the very special case $\vartheta
= k^0_n\Delta_n$, using $\widehat{k^0_n}$ defined as one maximizer in
$k$ of the contrast sequence
\[
k \leadsto|{\mathcal C}_n(k)|.
\]
Indeed, from the preceding computations, we have
\begin{equation}\label{almostsure}
\mbox{Almost surely, for large enough }n,\qquad\widehat{k^0_n}=k^0_n.
\end{equation}
This is the essence of our method. For an arbitrary $\vartheta$, we can
anticipate that an approximation of~$\vartheta$ taking the form
$k^0_n\Delta_n$ would add an extra error term of the order of the
approximation, that is, $\Delta_n$, which is a first guess for an
achievable rate of convergence.\vadjust{\goodbreak}

In a general context of regular semi-martingales with lead-lag
effect, sampled at random non-synchronous data points, we consider
the Hayashi--Yoshida (later abbreviated by HY) covariation estimator
and modify it with an appropriate time shift on one component. We
maximize the resulting empirical covariation estimator with respect
to the time shift over an appropriate grid.
%Our estimator
%is robust to nonsynchronous data and does not require any
%pre-processing of the data such as previous tick algorithm or so.
%The resulting
%and which is robust to
\subsubsection*{Construction of the estimator}
We need some notation. If $H = (\underline{H}, \overline{H}]$ is an
interval, for $\vartheta\in\Theta$, we define the shift interval
$H_\vartheta:=H+\vartheta= (\underline{H}+\vartheta, \overline
{H}+\vartheta]$. We write
\[
X(H)_t:=\int_{0}^t 1_H(s)\,\mathrm{d}X_s
\]
for a (possibly random) interval, such that $s \leadsto
1_{H}(s)$ is an elementary predictable process. Also, for notational
simplicity, we will often use the abbreviation
\[
X(H):=X(H)_{T+\delta} = \int_{0}^{T+\delta} 1_H(s)\,\mathrm{d}X_s.
\]
The shifted HY covariation contrast is defined as the function
\begin{eqnarray*}
\tilde\vartheta&\leadsto&{\mathcal U}^n(\tilde{\vartheta
})\\
&:=&1_{\tilde
{\vartheta}\geq0}\sum_{I \in
{\mathcal I}, J \in{\mathcal J}, \overline{I}\leq T}
X(I)Y(J)1_{\{I \cap J_{-\tilde{\vartheta}}\neq
\varnothing\}}\\
&&{}+1_{\tilde{\vartheta}<0}\sum_{I \in{\mathcal I}, J
\in
{\mathcal J}, \overline{J}\leq T} X(I)Y(J)1_{\{J \cap
I_{\tilde{\vartheta}}\neq\varnothing\}}.
\end{eqnarray*}
Our estimator $\widehat\vartheta_n$ is obtained by maximizing the
contrast $\tilde{\vartheta}\leadsto|{\mathcal U}^n(\tilde
{\vartheta
})|$ over the finite grid ${\mathcal G}^n$ constructed in
Assumption B3 in Section~\ref{statmodel} above. Eventually,
$\widehat
\vartheta_n$ is defined as a solution of
\begin{equation} \label{defestimator}
| {\mathcal U}^n(\widehat\vartheta_n)|=\max_{\tilde
{\vartheta
} \in{\mathcal G}^n}| {\mathcal U}^n(\tilde{\vartheta})|.
\end{equation}
%
%{\tt{\small[EXISTENCE AND UNIQUENESS TO BE DISCUSSED HERE]}}.
%s3.3 ###
\subsection{Convergence results} \label{th}
Since $\tau_{-\vartheta}(Y^c)$ is a $\mathbb{F}$-local martingale, the
quadratic variation process $\langle X^c, \tau_{-\vartheta
}(Y^c)\rangle
$ is well defined. We are now ready to assess our main result:
\begin{thmm} \label{mainth}
Work under Assumptions \textup{\ref{assa}} and \textup{\ref{assb}}. The estimator $\widehat\vartheta_n$
defined in \eqref{defestimator} satisfies
\[
v_n^{-1}(\widehat\vartheta_n - \vartheta) \rightarrow0
\]
in probability, on the event $\{\langle X^c,\tau_{-\vartheta
}(Y^c)\rangle_T \neq0\}$, as $n \rightarrow\infty$.
%Probleme avec la var quadratique
\end{thmm}

Theorem~\ref{mainth} provides a rate of convergence for our estimator:
the accuracy $\Delta_n^{-1}$ is nearly achievable, to within arbitrary
accuracy. The next logical step is the availability of a central limit
theorem. In the general case, this is not straightforward. We may,
however, be more accurate if we further restrict ourselves to
synchronous data in the Bachelier case; that is,
we have data
\begin{equation} \label{synchronous}
(X_0,Y_0), (X_{\Delta_n}, Y_{\Delta_n}), (X_{2\Delta_n}, Y_{2\Delta
_n}),\ldots
\end{equation}
over the time interval $[0,T]$, and the process $(X,Y)$ admits
representation \eqref{defquickbachelier}. We can then exhibit the
asymptotic behavior of the contrast function $\vartheta\leadsto
{\mathcal U}^n(\vartheta)$, in a vicinity of size $\Delta_n$, of the
lead-lag parameter.
More precisely, we have the following proposition.
%As a by-product, the heuristics of the preliminary of Section

\begin{propm} \label{asympt_behaviour_contrast}
Let $\varphi(t) = (1-|t|)1_{|t| \leq1}$ denote the usual hat function.
Let us consider the Bachelier model \eqref{defquickbachelier} and a
synchonous observation sampling scheme \eqref{synchronous}, with
lead-lag parameter $\vartheta\in\Theta$. If $|\tilde{\vartheta} -
\vartheta| \leq\Delta_n$, we have
\[
{\mathcal U}^n(\tilde{\vartheta}) = \sigma_1\sigma_2 \bigl(T\rho
\varphi\bigl(\Delta_n^{-1}(\tilde{\vartheta}-\vartheta)
\bigr)+T^{1/2}\Delta_n^{1/2}\sqrt{1+\rho^2\varphi\bigl(\Delta
_n^{-1}(\tilde
{\vartheta}-\vartheta)\bigr)} \xi^n\bigr),
\]
where $\xi^n$ is a sequence of random variables that converge in
distribution to the standard Gaussian law ${\mathcal N}(0,1)$ as $n
\rightarrow\infty$.
\end{propm}

This representation is useful to understand the behavior of the
contrast function ${\mathcal U}_n(\tilde{\vartheta})$: up to a scaling
factor, $|{\mathcal U}^n(\tilde{\vartheta})|$ is asymptotically
proportional to the realization of the absolute value\vspace*{1pt} of Gaussian
random variable $|{\mathcal N}(m_n(\tilde{\vartheta
}),a_n(\tilde
{\vartheta})^2)|$, with
\[
m_n(\tilde{\vartheta})=T\rho  \varphi\bigl(\Delta_n^{-1}(\tilde
{\vartheta}-\vartheta)\bigr) \quad \mbox{and}\quad
a_n(\tilde{\vartheta})=T^{1/2}\Delta_n^{1/2}\sqrt{1+\rho^2\varphi
\bigl(\Delta_n^{-1}(\tilde{\vartheta}-\vartheta)\bigr)}
\]
which has asymptotic value $m_n(\tilde{\vartheta})$ as soon as the mean
dominates the standard deviation. We then have
\[
\frac{|m_n(\tilde{\vartheta})|}{a_n(\tilde{\vartheta})} = \Delta
_n^{-1/2} \rho T^{1/2}\frac{\varphi(\Delta_n^{-1}(\tilde
{\vartheta
}-\vartheta))}{\sqrt{1+\rho^2 \varphi(\Delta
_n^{-1}(\tilde
{\vartheta}-\vartheta))}} \rightarrow\infty\qquad  \mbox{as }  n
\rightarrow\infty,
\]
and this is the case if $|\tilde{\vartheta}-\vartheta| \leq\Delta_n$;
otherwise, the pike $ \rho  \varphi(\Delta_n^{-1}(\tilde
{\vartheta
}-\vartheta))$ degenerates toward~$0$, and the contrast behaves
like a non-informative $\Delta_n^{1/2} |{\mathcal N}(0,1)|$ up
to a multiplicative constant. It is noteworthy that Proposition \ref
{asympt_behaviour_contrast} reveals the influence of the correlation
$\rho$ in the estimation procedure. We see that if $\rho$ is too small,
namely of order $\Delta_n^{1/2}$, the same kind of degeneracy
phenomenon occurs: we do not have the divergence $m_n(\tilde{\vartheta
})/a_n(\tilde{\vartheta}) \rightarrow\infty$ anymore, and both mean
and standard deviation are of the same order; in that latter case,
maximizing $|{\mathcal U}^n(\tilde{\vartheta})|$ does not locate the
true value $\vartheta$.

The situation is a bit more involved when looking further for the next
logical step, that is, a limit theorem for $\widehat\vartheta_n \in
\operatorname{argmax}_{\tilde{\vartheta} \in{\mathcal G}^n}|{\mathcal
U}^n(\tilde{\vartheta})|$. The function $\tilde{\vartheta}
\leadsto
{\mathcal U}^n(\tilde{\vartheta})$ is not smooth, even asymptotically:
up to normalizing by $\Delta_n^{-1}$, $\tilde{\vartheta} \leadsto
\varphi(\Delta_n^{-1}(\tilde{\vartheta}-\vartheta))$ weakly
converges to a Dirac mass at point $\vartheta$, see Proposition \ref
{asympt_behaviour_contrast}.
In that case, it becomes impossible, in general, to derive a simple
central limit\vadjust{\goodbreak} theorem for $\widehat\vartheta_n$. Consider again the
synchronous case over $[0,T]=[0,1]$, and pick a regular grid $\mathcal
{G}^n$ with mesh $h_n$ such that $h_n \Delta_n^{-1}$ goes to zero. In
this situation, the contrast function is constant over all the points
belonging to one given interval of the form $(i\Delta_n,(i+1)\Delta
_n)$, for $i\in\mathbb{Z}$. For definiteness and without loss of
generality, we set
\[
\widehat\vartheta_n=\operatorname{min}\Bigl\{\vartheta_n,~\vartheta_n \in
\mathop{\operatorname{argmax}}_{\tilde{\vartheta}\in\mathcal{G}^n}|\mathcal
{U}^n(\tilde{\vartheta})|\Bigr\}.
\]
From Theorem~\ref{mainth}, we know that $v_n^{-1}(\widehat\vartheta
_n-\vartheta)$ goes to zero for any sequence $v_n$ such that $v_n^{-1}
\Delta_n\rightarrow0$; therefore, we look for the behavior of the
normalized error, with rate $\Delta_n^{-1}$. However, the following
negative result shows that this cannot happen.
\begin{propm}\label{noclt}
Under the preceding assumptions, there is no random variable $Z$ such that
$\Delta_n^{-1}(\widehat\vartheta_n-\vartheta)$ converges in
distribution to $Z$.
\end{propm}

The proof is given in the \hyperref[appendix]{Appendix}. Proposition~\ref{noclt} stems from
the fact that part of the error of~$\widehat\vartheta_n$ is given by
the difference between $\vartheta$ and its approximation on the grid
$\mathcal{G}^n$. This error is deterministic and cannot be controlled
at the accuracy level $\Delta_n$; see the proof in the \hyperref[appendix]{Appendix}. This
phenomenon is somehow illustrated in the simulation in Section \ref
{illustration}.
Note that this negative result is not in contradiction to result \eqref
{almostsure} which states that almost surely, for large enough $n$,
\mbox{$\widehat\vartheta_n=\vartheta$}. Indeed, result \eqref{almostsure} is
obtained considering a grid with mesh $\Delta_n$ and a very special
sequence of models where~$\vartheta$ is of the form $\vartheta=\vartheta_n=k^0_n\Delta_n$, with
$k^0_n$ an integer. In the case where~$\vartheta$ does not depend on
$n$, one can, of course, extend the almost sure result \eqref
{almostsure}. However, what can be obtained is essentially that almost
surely, for large enough $n$, $\vartheta\in(\widehat\vartheta
_n-\Delta
_n,\widehat\vartheta_n+\Delta_n)$. Therefore, we almost surely identify
the interval of size $2\Delta_n$ in which $\vartheta$ lies, but our
method does not enable us to say something more accurate.
\subsection{Discussion} \label{discussion}

\subsubsection*{Covariation estimation of non-synchronous data}
The estimation of the covariation between two semi-martingales from
discrete data from non-synchronous observation times has some
history. It was first introduced by Hayashi and Yoshida~\cite{HY1}
and subsequently studied in various related contexts by several
authors. A comprehensive list of references include: Malliavin and
Mancino~\cite{MM}, Hayashi and Yoshida~\cite{HY1,HY2,HY3,HY4,HY5}, Hayashi and Kusuoka~\cite{HK},
Ubukata and Oya~\cite{UO}, Hoshikawa \textit{et al.}~\cite{Hetal} and
Dalalyan and Yoshida~\cite{DY}.

\subsubsection*{About the rate of convergence}
The condition $\Delta_n=\mathrm{o}(v_n)$ of Assumption B1 is needed for
technical reasons, in order to manage the fact that $\Delta_n$ is
random in general. In the case of regular sampling $\Delta_n= n^{-1}$
with $T=1$,
the nearly obtained rate $\Delta_n=n^{-1}$ is substantially better than
the usual
$n^{-1/2}$-rate of a regular parametric statistical model. This is
due to the fact that the estimation of the lead-lag parameter is
rather a change-point detection problem; see~\cite{IH}\vadjust{\goodbreak} for a
general reference for the structure of parametric models. A more
detailed analysis of the contrast function shows that its limit is
not regular (not differentiable in the $\vartheta$-variable), and
this explains the presence of the rate $n^{-1}$. However, the
optimality of our procedure is not granted, and the rate $\Delta_n$
could presumably be improved in certain special situations.\vspace*{-2pt}

\subsubsection*{Lead-lag effect and arbitrage}\vspace*{-1pt}
As stated, the lead-lag model for the two-dimensional process
$(X,Y)$ is not a semi-martingale, unless one component is
appropriately shifted in time. This is \textit{not compatible} in
principle with the dominant theory of no-arbitrage models. This kind
of modeling, however, seems to have some relevance in practice, and
there is a natural way to reconcile both points of view.

We focus, for example, on the simplest Bachelier model of Section
\ref{leadlag Bachelier}. We show in this paper that the lead-lag
parameter $\vartheta$ can almost be identified in principle.
Consequently, the knowledge of $\vartheta$ can then be
incorporated into a trading strategy. If $\vartheta\neq0$, we
can obtain, in principle, some \textit{statistical arbitrage}, in the
sense that we can find, in the Bachelier model without drift, a
self financing portfolio of assets $X$ and $\tau_{-\vartheta}(Y)$
with initial value zero and whose expectation at time~$T$ is
positive. %See also~\cite{G} for a general formulation of arbitrage
%in this context.

This statistical arbitrage can be erased by introducing further trading
constraints such as a maximal trading frequency and transaction cost
(slippage, execution risk and so on). In this setting, we can no longer
guarantee a statistical arbitrage. Moreover, we may certainly
incorporate risk constraints in order to define an admissible strategy.

This outlines that although we perturb the semi-martingale classical
approach, our lead-lag model is compatible in principle with
non-statistical arbitrage constraints, under refined studies of risk
profiles. We intend to set out, in detail, these possibilities in a
forthcoming work.\vspace*{-2pt}
%{\tt[INSERT MORE MATERIAL HERE]}

\subsubsection*{Microstructure noise}\vspace*{-1pt}
Our model does not incorporate microstructure noise. This is
reasonable if $\Delta_n$ is thought of on a daily basis, say (if $T$
is of the order of a year or more say), but is inconsistent in a
high-frequency setting where $T$ is of the order of one day. In that
context, efficient semi-martingale prices of the assets are subject
to the so-called microstructure noise; see, among others,
Zhang \textit{et al.}~\cite{ZMAb}, Bandi and Russell~\cite{BR},
Barndorff-Nielsen \textit{et al.}~\cite{BNHLS}, Hansen and Lunde
\cite{HL}, Jacod \textit{et al.}~\cite{Jetal}, Rosenbaum~\cite{R2,R}. In~\cite{RR1} and~\cite{RR2}, Robert and Rosenbaum
introduce a model (model with uncertainty zones) where the efficient
semi-martingale prices of the assets can be estimated at some random
times from the observed prices. In particular, it is proved that the
usual Hayashi--Yoshida estimator is consistent in this microstructure
noise context as soon as it is computed using the estimated values
of the efficient prices. Using the same approach, that is, applying
the lead-lag estimator to the estimated values of the efficient
prices, one can presumably build an estimator which is robust to
microstructure noise.\vspace*{-2pt}

\subsubsection*{How to use high-frequency data in practice}\vspace*{-1pt}

Nevertheless, when high-frequency data are considered, we propose a
simple pragmatic methodology that allows us to implement our lead-lag\vadjust{\goodbreak}
estimation procedure without requiring the \mbox{relatively} involved data
pre-processing suggested in the previous paragraph. A preliminary
inspection of the signature plot in trading time -- the realized
volatility computed with different subsampling values for the
trading times -- enables us to select a coarse subgrid among the
trading times where microstructure noise effects can be neglected.
Thanks to the non-synchronous character of high-frequency data, we
can take advantage of this subsampling in trading time and obtain
accurate estimation of the lead-lag parameter, at a scale that is
significantly smaller than the average mesh size of the coarse grid
itself. This would not be possible with a regular subsampling in
calendar, time where the price at time $t$ would be defined as the
last traded price before~$t$. This empirical approach is developed
in the numerical illustration Section~\ref{illustration real data}
on real data, in the particular case of measuring lead-lag between
the future contract on Dax (FDAX) and the Euro-Bund future contract
(FGBL) with same maturities.
%In this context, our approach shall be modified, and a
%pre-processing of the data is presumably required.
%An interesting possibility would be to incorporate pre-averaging
%in our data sets, as introduced by Jacod \textit{et al.}
%bounds.
\subsubsection*{Extension of the model}
We consider this work as a first -- and relatively simple -- attempt
for modeling the lead-lag effect in continuous time models. As a
natural extension, it would presumably be more reasonable to
consider more intricate correlations between assets in the model.
For example, one could add a common factor in the two assets,
without lead-lag effect, as suggested by the empirical study of
Section~\ref{illustration real data}. Through this, and in addition to
the ``lead-lagged correlation,'' one would also obtain an
instantaneous correlation between the assets. In order to estimate
the lead-lag parameter in this context, one would presumably be
required to consider local maxima of the contrast function we
develop here. Such a development is again left out for future work.

%s4 ###
\section{\texorpdfstring{Proof of Theorem \protect\ref{mainth}}{Proof of Theorem 1}}\label{proofs}

The proof of Theorem~\ref{mainth} is split in four parts. In the
first three parts, we work under supplementary assumptions on the
processes and the parameter space (Assumption \hyperref[assaa]{$\widetilde{\mathrm{A}}$}). We
first show that if we compute the contrast function over points
$\vartheta_n$ of the grid ${\mathcal G}^n$ such that the order of
magnitude of $|\vartheta_n-\vartheta|$ is bigger than $v_n$, then
the contrast function goes to zero (Proposition~\ref{step1}). Then
we prove that, on the contrary, if the order of magnitude of
$|\vartheta_n-\vartheta|$ is essentially smaller than $v_n$, then
the contrast function goes to the covariation between $X$ and
$\tau_{-\vartheta}(Y)$ (Proposition~\ref{step2}). We put these two
results together in the third part which ends the proof of Theorem
\ref{mainth} under the supplementary assumptions. The proof under
the initial assumptions is given in the last part.
%s4.1 ###
\subsection{Preliminaries}
\subsubsection*{Supplementary assumptions}
For technical convenience, we will first prove Theorem~\ref{mainth}
when the sign of $\vartheta$ is known and when the components $X$ and
$Y$ are local martingales. Moreover, we introduce a localization
tool.\vadjust{\goodbreak}
The quadratic variation processes of $X$ and $Y$ admitting locally
bounded derivatives, there exists a sequence of stopping times tending
almost surely to $T+\delta$ such that the associated stopped processes
are bounded by deterministic constants. Since Theorem~\ref{mainth} is a
convergence in probability result, we can, without loss of generality,
work under the supplementary assumption that the quadratic variation
processes are bounded over $[0,T+\delta]$. Therefore, we add-up the
following restrictions:

\renewcommand{\theass}{$\widetilde{\bolds{\Alph{ass}}}$}
\setcounter{ass}{0}
\begin{ass}\label{assaa}
We
have Assumption \textup{\ref{assa}} and:
\begin{enumerate}[$\widetilde{\mathrm{A}}1$.]
\item[$\widetilde{\mathrm{A}}1$.] There
exists $L>0$ such that $\langle X\rangle'_{T+\delta}\leq L$ and
$\langle
Y\rangle'_{T+\delta}\leq L$.
\item[$\widetilde{\mathrm{A}}2$.] The
parameter set is restricted to $\vartheta= [0,\delta)$.
Consequently, by $\mathcal{G}^n$ we mean here
$\mathcal{G}^n\cap[0,\delta)$.
\item[$\widetilde{\mathrm{A}}3$.] $X = X^c$ and
$Y=Y^c$.
\end{enumerate}
\end{ass}

\begin{notation*}
We now introduce further notation. For $I\in{\mathcal I}$ and $J \in
{\mathcal J}$, let
\[
\underline{I^n}=\underline{I}\wedge
\inf\Bigl\{t,
\max_{I'}\{\overline{I'}\wedge t-\underline{I'} \wedge t\}
\geq v_n \Bigr\}
\wedge T
\]
and %Similarly, for an $\F^{\vartheta^*}$-stopping time $\tau_n$,
\[
\underline{J^n}=\underline{J}\wedge
\inf\Bigl\{t,
\max_{J'}\{\overline{J'}\wedge t-\underline{J'} \wedge t\}
\geq v_n \Bigr\}
\wedge(T+\delta).
\]
We define $\overline{I^n}$ and $\overline{J^n}$ in the same way for
$\overline{I}$ and $\overline{J}$, respectively.
Let $I^n=(\underline{I^n},\overline{I^n}]$ and
$J^n=(\underline{J^n},\overline{J^n}]$.
\end{notation*}
\begin{rem}
We have the following interpretation of
$\underline{I^n}$ and $\overline{I^n}$: let
$\tau^n$ denote the first time for which we know that an interval $I$
will have a width that is larger than $v_n$. Then we keep only the
$\underline{I}$ and $\overline{I}$ that are smaller than $\tau^n$. If
$\tau^n\leq T$, we also consider $\tau^n$
among the observation times. Note that $\tau^n$ is not a true
observation time in general. However, this will not be a problem
since the set where $\Delta_n$ is bigger than $v_n$ will be
asymptotically negligible. Obviously $\underline{I^n}$ and
$\overline{I^n}$ are $\F$-stopping times, and $\underline{J^n}$
and $\overline{J^n}$ are $\mathbb{F}^{\vartheta+v_n}$-stopping
times.
\end{rem}

Finally, for two intervals $H = (\underline{H},
\overline{H}]$ and $H'=(\underline{H'}, \overline{H'}]$, we define
\[
K(H,H') := 1_{H\cap H' \neq\varnothing}.
\]

%s4.2 ###
\subsection{The contrast function}

We consider here the case where the order of magnitude of
$|\vartheta_n-\vartheta|$ is bigger than $v_n$. We first need to give a
preliminary lemma that will ensure
that the quantities we will use in the following are well defined.

\begin{lemm}\label{lagest-1} Work under Assumption \textup{B2}, under the
slightly more general assumption that for all $I =(\underline{I},
\overline{I}] \in{\mathcal I}$, the random variables $\underline{I}$
and $\overline{I}$ are $\F$-stopping times. Suppose that $\widetilde
\vartheta\geq\vartheta+\ep_n$
and $2v_n\leq\ep_n$. Then for any random variable $X'$\vadjust{\goodbreak} measurable
w.r.t. $\calf_{\overline{I^n}}$, the random variable
$X'K(I^n_{\widetilde\vartheta},J^n)$ is
$\calf^{\vartheta}_{\underline{J^n}}$-measurable. In particular,
$f(\overline{I^n})X(I^n)K(I^n_{\widetilde\vartheta},J^n)$ is
$\calf^{\vartheta}_{\underline{J^n}}$-measurable
for any measurable function $f$.
\end{lemm}

The proof of Lemma~\ref{lagest-1} is given in the \hyperref[appendix]{Appendix}. It is
important to note that Lemma~\ref{lagest-1} implies that for
$\widetilde\vartheta\geq\vartheta+\ep_n$ and $2v_n\leq\ep_n$,
the random variable
\[
1_{\{\overline{I^n}\leq T\}}X(I^n)
K(I^n_{\widetilde\vartheta},J^n)1_{J^n}(s)
\]
is
$\calf^{\vartheta}_s$-measurable. Indeed, $1_{J^n}(s)$ is
$\calf^{\vartheta}_s$ and $1_{J^n}(s)=1$ implies
$s\geq\underline{J^n}$. We now introduce a functional version of
${\mathcal U}^n$ by considering the random process
\[
\bbU^n(\widetilde\vartheta)_t:=\sum_{I\in{\mathcal I}, J\in
{\mathcal
J}, \overline{I^n}\leq T}X(I^n)Y(J^n)_tK(I_{\widetilde\vartheta}^n, J^n).
\]
We are now able to give the main proposition for the vanishing of
the contrast function.

\begin{propm}\label{step1}
Let $\ep_n=2v_n$, $ \calg^n_+=\{\widetilde\vartheta\in\calg^n,
\widetilde\vartheta\geq\vartheta+\ep_n\}$ and
$\calg^n_-=\{\widetilde\vartheta\in\calg^n, \widetilde
\vartheta\leq\vartheta-\ep_n\}$. We have
\[
\max_{\widetilde\vartheta\in\calg^n_+\cup\calg^n_-}|\bbU
^n(\widetilde\vartheta)_{T+\delta}|\rightarrow0,
\]
in probability.
\end{propm}
\begin{pf}
Assume first $\widetilde\vartheta\geq\vartheta+\ep_n$. Thanks to Lemma
\ref{lagest-1}, we obtain a martingale representation of the
process $\bbU^n(\widetilde\vartheta)$ that takes the form

\[
\bbU^n(\widetilde\vartheta)_t=\sum_{I \in{\mathcal I},J \in
{\mathcal J}}\int_0^t1_{\{\overline{I^n}\leq T\}} X(I^n)
K(I^n_{\widetilde\vartheta},J^n) 1_{J^n}(s)\,\mathrm{d} Y_s,
\]
where the
stochastic integral with respect to $Y$ is taken for the
filtration $\F^{\vartheta}$.
As a result, the $\F^{\vartheta}$-quadratic variation of $\bbU^n$ is
given by
\[
\langle\bbU^n(\widetilde\vartheta) \rangle_t = \int_0^t
\biggl(\sum_{I \in{\mathcal I},J \in{\mathcal J}}
1_{\{\overline{I^n}\leq T\}} X(I^n)K(I^n_{\widetilde
\vartheta},J^n)1_{J^n}(s)\biggr)^2\,\mathrm{d}\langle Y\rangle_s.
\]
Using that the
intervals $J^n$ are disjoint, we obtain
\[
\langle\bbU^n(\widetilde\vartheta) \rangle_t = \int_0^t \sum_{J
\in{\mathcal J}}\biggl(\sum_{I \in{\mathcal I}}
1_{\{\overline{I^n}\leq T\}} X(I^n)K(I^n_{\widetilde
\vartheta},J^n)\biggr)^21_{J^n}(s)\,\mathrm{d}\langle Y\rangle_s.
\]
For a given
interval $J^n$, the union of the intervals $I^n$ that have a non-empty
intersection with $J^n$ is an interval of width smaller than
$3 v_n$. Indeed, the maximum width of $J^n$ is $v_n$ and add to
this (if it exists) the width of the interval $I^n$\vadjust{\goodbreak} such that
$\underline{I^n}\leq\overline{J^n}$,
$\overline{I^n}\geq\overline{J^n}$ and the width of the interval~$I^n$ such that $\underline{I^n}\leq\underline{J^n}$,
$\overline{I^n}\geq\underline{J^n}$. Thus,
\begin{eqnarray*}
\sum_{I \in{\mathcal I}}1_{\{\overline{I^n}\leq T\}
}X(I^n)K(I^n_{\widetilde
\vartheta},J^n)&\leq&
\mathop{\operatorname{sup}}_{s\leq T}\mathop{\operatorname{sup}}_{0\leq u\leq3
v_n}\bigl|X_{(s+u)\wedge T}-X_u\bigr|\\[-2pt]
&\leq&2\max_{1\leq k\leq\lfloor(3 v_n)^{-1}T\rfloor} \sup_{t\in[3
v_n(k-1),3 v_nk]}\bigl|X_{t\wedge T}-X_{3v_n(k-1)}\bigr|.
\end{eqnarray*}
Consequently, we obtain for every $t\in[0,T+\delta]$ and $\widetilde
\vartheta\in[\vartheta+\ep_n,\delta]$,
\[
\langle\bbU^n(\widetilde\vartheta) \rangle_t\leq4L
(T+\delta)\max_{1\leq k\leq\lfloor(3
v_n)^{-1}T\rfloor}\sup_{t\in[3 v_n(k-1),3 v_nk]} \bigl|X_{t\wedge
T}-X_{3(k-1)v_n}\bigr|^2.
\]
For every $p>1$, it follows from the
B\"urkholder--Davis--Gundy inequality that
\[
\E[|\bbU^n(\widetilde
\vartheta)_{T+\delta}|^{2p}]\lesssim\sum_{k=1}^{\lfloor(3
v_n)^{-1}T\rfloor}\E\Bigl[\sup_{t\in[3 v_n(k-1),3 v_nk]} \bigl|X_{t\wedge
T}-X_{3(k-1)v_n}\bigr|^{2p}\Bigr]\lesssim v_n^{ p-1},
\]
where the symbol
$\lesssim$ means inequality in order, up to constant that does not
depend on~$n$. Pick $\varepsilon
>0$. We derive
\begin{eqnarray*}
\PP\Bigl[\max_{\widetilde\vartheta
\in\calg^n_+}|\bbU^n(\widetilde\vartheta)_{T+\delta}|>\ep\Bigr]
&\leq&
\ep^{-2p}\sum_{\widetilde
\vartheta\in\calg^n_+}\E[|\bbU^n(\widetilde
\vartheta)_{T+\delta}|^{2p}]
\\[-2pt]
&\lesssim &v_n^{p-1} \sharp \calg^n_+ \to0
\end{eqnarray*}
as $n
\rightarrow\infty$, provided $p>\gamma+1$ where $\gamma$ is
defined in Assumption B3, a choice that is obviously possible. The
same argument holds for the case $\widetilde\vartheta\leq
\vartheta-\ep_n$, but with an $X$-integral representation in that
latter case. The result follows.
\end{pf}

%s4.3 ###
\subsection{Stability of the HY estimator}

We consider now the case where the order of magnitude of
$|\vartheta_n-\vartheta|$ is essentially smaller than~$v_n$. We
have the following proposition:

\begin{propm} \label{step2}
Work under Assumptions \textup{\hyperref[assaa]{$\widetilde{\mathrm{A}}$}} and \textup{\ref{assb}}. For any sequence
$\vartheta_n$ in $[0,\delta)$ such that $\vartheta_n\leq
\vartheta$ and $|\vartheta_n -\vartheta|\leq\rho_n$ (remember
that $\rho_n$ is defined in Assumption \textup{B3}), we have
\begin{equation} \label{convstability}
\mathcal{U}^n(\vartheta_n)\rightarrow\langle
X,\tau_{-\vartheta}(Y)\rangle_{[0,T]},
\end{equation}
in probability as $n \rightarrow\infty$.
\end{propm}

\begin{pf}
The proof goes into several steps.

\textit{Step 1.} In this step, we show that our contrast
function can be regarded as the Hayashi--Yoshida estimator applied to
$X$ and to the properly shifted values of $Y$ plus a remainder
term. If $\vartheta=0$, then $\vartheta_n=0$, and Proposition
\ref{step2} asserts nothing but the consistency of the standard
HY-estimator; see Hayashi and Yoshida~\cite{HY1} and Hayashi and
Kusuoka~\cite{HK}. Thus we may assume $\vartheta>0$.\vadjust{\goodbreak}

By symmetry, we only need to consider the case where
%$$\E[\sharp\{I\}] \geq
%$$
%
\[
\E[\sharp{\mathcal I}] \geq
\E[\sharp{\mathcal J}].
\]
Set $\delta_n=\vartheta_n-\vartheta$, $\widetilde Y_t
 =\tau_{-\vartheta}(Y)_t$ and $\circj=J^n_{-\vartheta}$ and
\[
\mathbb{U}^n(\vartheta_n)=\sum_{I \in{\mathcal I},J \in{\mathcal
J},\overline{I}\leq
T}X(I^n)Y(J^n) 1_{\{I^n\cap J^n_{-\vartheta_n}\not=\varnothing\}}.
\]
We then have
\[
\mathbb{U}^n(\vartheta_n)=\sum_{I \in{\mathcal I},J \in{\mathcal
J},\overline{I}\leq T}X(I^n)\widetilde Y(\circj) 1_{\{I^n\cap
\circj_{-\delta_n}\not=\varnothing\}}.
\]
This can be written ${\mathcal V}^n+{\mathcal R}^n$ with
\begin{eqnarray*}
{\mathcal
V}^n &=&\sum_{I \in{\mathcal I},J \in{\mathcal
J},\overline{I}\leq T}X(I^n)\circy(\circj_{-\delta_n})
1_{\{I^n\cap\circj_{-\delta_n}\not=\varnothing\}},\\[-2pt]
{\mathcal R}^n &=& \sum_{I \in{\mathcal I},J \in{\mathcal
J},\overline{I}\leq T}X(I^n)
\{\circy(\circj)-\circy(\circj_{-\delta_n})\} 1_{\{
I^n\cap
\circj_{-\delta_n}\not=\varnothing\}}.
\end{eqnarray*}

Remark that $\circy(\circj)$ and $\circy(\circj_{-\delta_n})$ are well
defined since $\circy$ is defined on $[-\vartheta,T]$ and $\vartheta
_n\leq\vartheta$.
For every $J \in{\mathcal J}$, $\overline{J^n}$ is a
$\F^{\vartheta+v_n}$-stopping time; therefore
$\overline{\circj_{-\delta_n}}=\overline{J^n_{-\vartheta_n}}$ is a
$\F^{v_n-\delta_n}$-stopping time, and a $\F$-stopping time as
well.
%It is in force for $\overline{\circj}$.
Thus ${\mathcal V}^n$ is a variant of the HY-estimator: more
precisely,
\[
\widetilde{\mathcal V}^n :=
\sum_{I,J:\overline{I}\leq
T}X(I^n)\circy(\circj_{-\delta_n}\cap \bbR_+) 1_{\{I^n\cap
(\circj_{-\delta_n}\cap \bbR_+)\not=\varnothing\}}
\]
is the
original HY-estimator, and we have ${\mathcal V}^n-\widetilde
{\mathcal V}^n \rightarrow0$ in probability as $n \rightarrow
\infty$. It follows that $\widetilde{\mathcal V}^n \rightarrow
\langle X,\circy\rangle_T$ in probability as $n\to\infty$; see
\cite{HY1} and~\cite{HK}.

\textit{Step 2.} Before turning to the term ${\mathcal
R}^n$, we give a technical lemma and explain a simplifying
procedure. For an interval $I = [\underline{I}, \overline{I}) \in
{\mathcal I}$, set
\[
M^I := \sup\{\overline{\circj_{-\delta_n}}, J \in{\mathcal
J}, \overline{\circj_{-\delta_n}}\leq\overline{I^n}\}.
\]
%
%When there is no $J$ satisfying
%$\overline{\circj_{-\delta_n}}\leq\overline{I^n}$,
%we set $L^I=0$.
Note that if we consider the interval $J$ at the extreme left end
of the family ${\mathcal J}$, we have, for large enough $n$,
\[
\overline{\circj_{-\delta_n}}\leq v_n-\vartheta-\delta_n
<-\frac{\vartheta}{2},
\]
say, so we may assume that
%an interval $J$
%satisfying $\overline{\circj_{-\delta_n}}\leq\overline{I^n}$
the set over which we take the supremum is non-empty.
%always exists.

%The notation of Section~\ref{} is in force here.
%
\begin{lemm} \label{stoppingtime2}
Work under Assumption \textup{B2}. The random variables $M^I$ are
%For $I \in{\mathcal I}$, let
%$$
%M^I =
%$$
%When there is no $J$ satisfying
%$\overline{\circj_{-\delta_n}}\leq\overline{I^n}$,
%we set $L^I=0$.
%$M^I$ are
%$\F'$-stopping times for $I \in{\mathcal I}$.
${\mathbb F}$-stopping times.
\end{lemm}

The proof of this lemma is given in the \hyperref[appendix]{Appendix}. We now use a
simplifying operation. For each $I^n$, we merge all the $J^n$ such
that $\circj_{-\delta_n} \subset I^n$. We call this procedure
$\Pi$-reduction. The $\Pi$-reduction produces a new sequence of
increasing random intervals extracted from the original sequence
$(\circj_{-\delta_n})$, which are $\F$-predictable by Lemma
\ref{stoppingtime2}. More precisely, the end-points are
$\F$-stopping times. It is important to remark that the
$\Pi$-reduction implies that there are at most two points of type
$J$ between any $\underline{I^n}$ and $\overline{I^n}$. Moreover,
since ${\mathcal R}^n$ is a bilinear\vspace*{1pt} form of the increments of $X$
and $\circy$, it is invariant under $\Pi$-reduction. Likewise for
the maximum length $\Delta_n$. Thus, without loss of generality,
we may assume that the $\circj_{-\delta_n}$ are $\Pi$-reduced.

\textit{Step 3.} We now turn to ${\mathcal R}^n$. We write
\[
I^n(\circj_{-\delta_n})=\bigcup_{{I \in{\mathcal I},\overline
{I}\leq
T,\{I^n\cap\circj_{-\delta_n}\not=\varnothing\}}} I^n.
\]
We have
\[
|{\mathcal R}^n|\leq\sum_{J\in\mathcal{J}}|\circy(\circj)-\circy
(\circj
_{-\delta_n})||X(I^n(\circj_{-\delta_n}))|.
\]
We now index the intervals $\widetilde J^n$ by $j$ and set
$\widetilde J^n_j=\{0\}$ if $j>\sharp\{J\}$. Thus, the preceding
line can be written
\[
|{\mathcal R}^n|\leq\sum_{j}|\circy(\circj_j)-\circy(\circj
_{-\delta
_n,j})||X(I^n(\circj_{-\delta_n,j}))|.
\]
Then the Cauchy--Schwarz inequality gives that $(\E[|{\mathcal
R}^n|])^2$ is smaller than
\[
\sum_{j}\E[|\circy(\circj_j)-\circy(\circj_{-\delta
_n,j})|^2]\sum_{j}\E
[|X(I^n(\circj_{-\delta_n,j}))|^2].
\]
We easily get that
\[
\sum_{j}\E[|\circy(\circj_j)-\circy(\circj_{-\delta
_n,j})|^2]\lesssim
\delta_n\sharp{\mathcal J},
\]
and we claim that (see next step)
\begin{equation} \label{keypropineq}
\sum_{j}\E[|X(I^n(\circj_{-\delta_n,j}))|^2]\lesssim1.
\end{equation}
Since $\delta_n\leq\rho_n$, Proposition~\ref{step2} readily
follows.

\textit{Step 4.} It remains to prove \eqref{keypropineq}. Here
we extend $(X_t)_{t\in\bbR_+}$ as $X_s=0$ for $s<0$, and
denote the extended one by the same ``$X$.''
This extension is just
for notational convenience, and causes no problem because,
in what follows, we use the martingale property of $X$
only over the time interval $\bbR_+$. For ease of notation, we also
stop writing the index $j$ for the intervals.
We begin with the following remark. Take an interval
$\circj_{-\delta_n}$, say $(J_1,J_2]$ and
$I^n(\circj_{-\delta_n})$ associated, say $(I_1,I_2]$. Call $J_0$
the last observation point of type $J$ occurring before $J_1$
and
$J_{-1}$ the last observation point of\vadjust{\goodbreak}
type $J$ occurring before $J_0$. Two situations are possible:

\begin{itemize}[-]
 \item[-] If there is no observation point of type $I$ between
$I_1$ and $I_2$, then, if it exists, $J_0$ is necessarily before
$I_1$. If it does not exist, we have $J_1\leq v_n$.

 \item[-] If there are some observation points of type $I$
between $I_1$ and $I_2$, then $J_0$ might also be between $I_1$
and $I_2$. However, thanks to the $\Pi$-reduction, we know that
$J_{-1}$ is necessarily smaller than $I_1$. Consequently, we have that
$|X(I^n(\circj_{-\delta_n}))|$ is smaller than\looseness=-1
\[
\sup_{t\in[\underline{{{\widetilde J^{n,-1}}_{-\delta
_n}}},\overline{I^{n,-1}]}}
|X_t-X_{\underline{{{\widetilde J^{n,-1}}_{-\delta_n}}}}|
+\sup_{t\in[\underline{{{\widetilde J^{n,-2}}_{-\delta
_n}}},\overline
{I^{n,-2}}]}
|X_t-X_{\underline{{{\widetilde J^{n,-2}}_{-\delta_n}}}}|+
\sup_{t\in[\underline{\circj_{-\delta_n}},\overline{I^n_+}]}
|X_t-X_{\underline{\circj_{-\delta_n}}}|,
\]\looseness=0
where we used the following notation:
\begin{itemize}[-]
\item[-] $I^n_+$ is the first interval $I^n$ such that
$\overline{I^n}$ exits to the right of $\circj_{-\delta_n}$.
\item[-] ${\widetilde J^{n,-1}_{-\delta_n}}$ denotes the interval of
the form $\circj_{-\delta_n}$ which is the nearest neighbor
to $\circj_{-\delta_n}$ on the left.
\item[-] ${\widetilde J^{n,-2}}_{-\delta_n}$ denotes the interval of
the form $\circj_{-\delta_n}$ which is the nearest neighbor
to ${\widetilde J^{n,-1}}_{-\delta_n}$ on the left.
\item[-] $I^{n,-1}$ is the first exit time to the right
of ${\widetilde J^{n,-1}}_{-\delta_n}$ among the $I^n$.
\item[-] $I^{n,-2}$ is the first exit time to the right
of ${\widetilde J^{n,-2}}_{-\delta_n}$ among the $I^n$.
\item[-] For $k=1,2$, if ${\widetilde J^{n,-k}}_{-\delta_n}$ is not
defined, $ \sup_{t\in[\underline{{{\widetilde
J^{n,-k}}_{-\delta_n}}},\overline{I^{n,-k}]}}
|X_t-X_{\underline{{{\widetilde J^{n,-k}}_{-\delta_n}}}}|=0.$
\end{itemize}
\end{itemize}
Hence we obtain
\[
\sum_j\E[|X(I^n(\circj_{-\delta_n,j}))|^2]
\lesssim\sum_j\E\Bigl[\sup_{t\in[\underline{\circj_{-\delta
_n}},\overline
{I^n_+}]}
|X_t-X_{\underline{\circj_{-\delta_n}}}|^2\Bigr]
\]
and so $ \sum_j\E[|X(I^n(\circj_{-\delta
_n,j})
)|^2]$ can be bound in order by
\[
\sum_j\E\Bigl[\Bigl(\sup_{t\in[\underline{\circj_{-\delta
_n}},\overline
{\circj_{-\delta_n}}]}
|X_t-X_{\underline{\circj_{-\delta_n}}}|\Bigr)^2\Bigr]+\E\biggl[
\biggl(\sum
_j\sup_{t\in[\underline{I^n_+},\overline{I^n_+}]}
|X_t-X_{\underline{I^n_+}}|\biggr)^2\biggr].
\]
Thanks to the $\Pi$-reduction, we know that a given interval of
the form $(\underline{I^n_+},\overline{I^n_+}]$ can be associated
to, at most, two values of type $J$. Thus the second term of the
preceding quantity is
smaller than
\[
2\E\biggl[\biggl(\sum_i\sup_{t\in[\underline{I^n},\overline{I^n}]}
|X_t-X_{\underline{I^n}}|1_{i\leq
\sharp\mathcal{I}}\biggr)^2\biggr],
\]
where $i$ is an indexing of the
intervals $[\underline{I^n},\overline{I^n})$. Note that each
$\underline{I^n_+}$ is an $\F$-stopping time as it is the maximum
among all
\[
\overline{I^n} \leq\overline{\circj_{-\delta_n}},
\]
together
with a strong predictability property; see Lemma
\ref{stoppingtime2} for a similar statement. So, using
B\"urkholder--Davis--Gundy inequality, \eqref{keypropineq} is proved
and Proposition~\ref{step2} follows.
\end{pf}

%s4.4 ###
\subsection{\texorpdfstring{Completion of proof of Theorem \protect\ref{mainth} under Assumption \protect\hyperref[assaa]{$\widetilde{\mathrm{A}}$}}
{Completion of proof of Theorem 1 under Assumption A}}

Write ${\mathcal A}=\{\langle X,\tau_{-\vartheta}(Y)\rangle_T
\not=0\}$. By Assumption B3, we have
\[
\bigcup_{\widetilde
\vartheta\in\calg^n}[\widetilde\vartheta-\rho_n,\widetilde
\vartheta+\rho_n)\supset\vartheta.
\]
Therefore, there exists a
sequence $\vartheta_n$ in $\calg^n$ such that $\vartheta_n \leq
\vartheta$ and $|\vartheta_n-\vartheta|\leq2\rho_n$. For
sufficiently large $n$, we have $\rho_n \leq\ep_n=2v_n$. Moreover,
on the event ${\mathcal A},$
\[
{\mathcal U}^n(\hat{\vartheta}_n)>\sup_{\widetilde\vartheta\in
\calg
^n_+\cup \calg^n_-}
| {\mathcal U}^n(\widetilde\vartheta)|
\]
implies
$|\hat{\vartheta}_n-\vartheta|<\ep_n$. It follows that
\[
\PP[\{|\hat{\vartheta}_n-\vartheta|\geq\ep_n\}\cap
{\mathcal A}]
\leq\PP\Bigl[ \Bigl\{ \sup_{\widetilde\vartheta\in\calg^n_+\cup
\calg^n_-}
| {\mathcal U}^n(\widetilde
\vartheta)|\geq| {\mathcal U}^n(\vartheta_n)|
\Bigr\}\cap{\mathcal A} \Bigr].
\]
Let $\ep>0$. For large enough
$n$, the probability to have $\Delta_n$ smaller than $v_n$ is
larger than $1-\ep$ and, consequently,
\[
\PP[\{|\hat{\vartheta}_n-\vartheta|\geq\ep_n\}\cap{\mathcal
A}]\leq\PP\Bigl[ \Bigl\{ \sup_{\widetilde
\vartheta\in\calg^n_+\cup \calg^n_-} |\bbU^n(\widetilde
\vartheta)_{T+\delta}|\geq|{\mathcal
U}^n(\vartheta_n)| \Bigr\}\cap{\mathcal A} \Bigr] +\ep.
\]
This
can be bounded in order by
\begin{eqnarray*}
\PP\biggl[|{\mathcal U}^n(\vartheta_n)|<\frac{1}{2}|\langle
X,\tau_{-\vartheta}(Y)\rangle_T|\biggr]
+ \PP\biggl[ \biggl\{
\sup_{\widetilde\vartheta\in\calg^n_+\cup \calg^n_-}
|\bbU^n(\widetilde\vartheta)_{T+\delta}|
>\frac{1}{2}|\langle X,\tau_{-\vartheta}(Y) \rangle_T|\biggr\}\cap
{\mathcal A} \biggr]+\ep,
\end{eqnarray*}
and this last quantity converges to $\varepsilon$ as $n\to\infty$
by applying Proposition~\ref{step1} and Proposition~\ref{step2}.

%s4.5 ###
\subsection{The case with drifts}

We now give the proof of Theorem~\ref{mainth} under Assumptions
$\widetilde {\mathrm{A}}1$, $\widetilde {\mathrm{A}}2$ and~\ref{assb}. The contrast ${\mathcal
U}^n(\widetilde\vartheta)$ admits the decomposition $ {\mathcal
U}^n(\widetilde\vartheta)
%&=&
%1_{\{I\cap J_{-\vartheta}\not=\varnothing\}}
= \widetilde{\mathcal U}^n(\widetilde\vartheta) + {\mathcal
R}^n(\widetilde\vartheta)
$
with
\[
\widetilde{\mathcal U}^n(\widetilde\vartheta) =
\sum_{I \in{\mathcal I}, J \in{\mathcal J},\overline{I}\leq T
}X^c(I)Y^c(J) 1_{\{I\cap J_{-\widetilde\vartheta}\not=\varnothing\}}
\]
and
\[
{\mathcal R}^n(\widetilde\vartheta) = \sum_{I
\in{\mathcal I},J \in{\mathcal J}} \bigl(X(I)B(J)+A(I)Y^c(J)\bigr)
1_{\{I\cap J_{-\widetilde\vartheta}\not=\varnothing\}}.
\]
For a function $t \rightarrow Z_t$ defined on the interval $H$,
introduce the modulus of continuity
\[
w_Z(a,H)=\sup\{|Z_t-Z_s|, s,t\in H,
|s-t|<a\},\qquad  a >0.\vadjust{\goodbreak}
\]
%
% for a function $Z$
%defined on the interval $H$.
We have
\[
\sup_{\widetilde\vartheta\in[0,\delta)} |{\mathcal
R}^n(\widetilde\vartheta)| \leq
 w_X(3\Delta_n, [0,T])\sup_{t\in[0,T+\delta]}|B_t|
+w_{Y^c}(3\Delta_n, [0,T+\delta])\sup_{t\in[0,T]}|A_t|,
\]
and this
term goes to $0$ in probability as $n\to\infty$.

Finally, the result is obtained in a similar way as in the no-drift
case, using $(X^c,Y^c)$ in place of $(X,Y)$.

%s4.6 ###
\subsection{\texorpdfstring{The case where $\vartheta\in(-\delta,\delta)$}
{The case where theta in (-delta,delta)}}

We now give the proof of Theorem~\ref{mainth} under Assumptions
$\widetilde {\mathrm{A}}1$ and~\ref{assb}. Even in the case where $\vartheta$ is
negative, Proposition $\ref{step1}$ is still in force, and we
obtain
\[
\mathop{\operatorname{sup}}_{\widetilde\vartheta\in{\mathcal
G}^n\cap[0,\delta)}|\mathcal{U}^n(\widetilde\vartheta
)|\rightarrow
0
\]
in probability as $n \rightarrow\infty$. The result follows from Remark
\ref{inverse}.
%{\tt[INSERT PROOF OF LOWER BOUND HERE]}

%s5 ###
\section{A numerical illustration on simulated data} \label{illustration}
%s5.1 ###
\subsection{Synchronous data: Methodology} \label{DGP}
We first superficially analyze the performances of $\widehat
\vartheta_n$ on a simulated lead-lag Bachelier model without
drift. More specifically, we take a random process $(X,
\tau_{-\vartheta}(Y))$ following the representation given in
\eqref{leadlagrepresentation} in Section~\ref{leadlag Bachelier},
having
\[
T=1, \qquad \delta= 1,\qquad  \vartheta= 0.1,\qquad x_0=\widetilde y_0=0,\qquad \sigma
_1=\sigma_2=1.
\]
In this simple model, we consider again synchronous, equispaced
data with period $\Delta_n$ and correlation parameters $\rho$. In
that very simple model, we construct $\widehat\vartheta_n$ with a
grid ${\mathcal G}^n$ with equidistant points with mesh\footnote{Note
that, strictly speaking, such grid is not fine enough in order to
fulfill our assumptions. However, the contrast function is constant
over all the points of a given
interval $(k\Delta_n,(k+1)\Delta_n)$, $k\in\mathbb{Z}$, and its
value is
just the sum of the values obtained for the shifts $k\Delta_n$ and
$(k+1)\Delta_n$.} $h_n=\Delta_n$. We consider the following variations:
\begin{enumerate}
\item Mesh size: $h_n \in\{10^{-3}, 3. 10^{-3}, 6. 10^{-3}\}$.
\item Correlation value: $\rho\in\{0.25, 0.5, 0.75\}$.
\end{enumerate}

%%
%t1 ###
\begin{table}
\def\arraystretch{0.9}
\caption{Estimation of $\vartheta=0.1$ on 300 simulated
samples for $\rho\in\{0.25, 0.5, 0.75\}$}\label{tab1}
\begin{tabular*}{\textwidth}{@{\extracolsep{\fill}}llllll@{}}
\hline
% after \  \hline or \cline{col1-col2} \cline{col3-col4} ...
$\widehat\vartheta_n$ &0.096 & 0.099 & 0.1& 0.102 & Other\\
\hline
FG, $\rho=0.75$ &\phantom{0}0& \phantom{00}0 & 300 & \phantom{00}0&\phantom{00}0\\
MG, $\rho=0.75$ &\phantom{0}0& 300 & \phantom{00}0 & \phantom{00}0 &\phantom{00}0\\
CG, $\rho=0.75$ &\phantom{0}1& \phantom{00}0 & \phantom{00}0 & 299 &\phantom{00}0\\[3pt]
FG, $\rho=0.50$ &\phantom{0}0& \phantom{00}0 & 300 & \phantom{00}0 &\phantom{00}0\\
MG, $\rho=0.50$ &\phantom{0}0& 299 & \phantom{00}0 & \phantom{00}1 &\phantom{00}0\\
CG, $\rho=0.50$ &13& \phantom{00}0 & \phantom{00}0 & 280 &\phantom{00}7\\[3pt]
FG, $\rho=0.25$ &\phantom{0}0& \phantom{00}0 & 300 & \phantom{00}0 &\phantom{00}0\\
MG, $\rho=0.25$ &\phantom{0}0& 152 & \phantom{00}0 & \phantom{0}11 &137\\
CG, $\rho=0.25$ &10& \phantom{00}0 & \phantom{00}0 & \phantom{0}66 &124\\
\hline
\end{tabular*}  \vspace*{-3pt}
\end{table}

%
%f1 ###
\begin{figure}[b]
\vspace*{-3pt}
\includegraphics{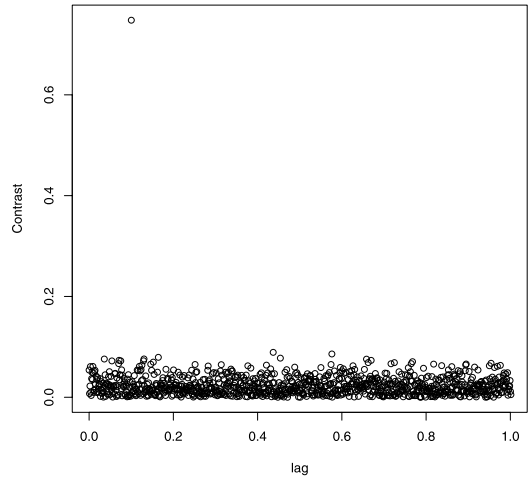}

\caption{Fine grid case (FG). Over one simulation:
displayed values of $| {\mathcal U}^n(\widetilde
\vartheta)|$ for $\widetilde\vartheta\in{\mathcal G}^n$
with mesh $h_n=10^{-3}$ and $\rho=0.75$. The value
$\max_{\widetilde\vartheta\in{\mathcal G}^n}| {\mathcal
U}^n(\widetilde\vartheta)|$ is well located.}\label{fig1}
\end{figure}

%s5.2 ###
\subsection{Synchronous data: Estimation results and their analysis}
We repeat 300 simulations of the experiment and compute the value
of $\widehat\vartheta_n$ each time, the true value being
$\vartheta=0.1$, letting $\rho$ vary in $\{0.25, 0.5, 0.75\}$. We
adopt the following terminology:\vadjust{\goodbreak}
\begin{enumerate}
\item The fine grid estimation (abbreviated FG) with $h_n=10^{-3}$.
\item The moderate grid estimation (abbreviated MG) with $h_n=3. 10^{-3}$.
\item The coarse grid estimation (abbreviated CG) with $h_n=6. 10^{-3}$.
\end{enumerate}
%
%mesh$=1/h=1/1000$, lead-lag$=0.1$, correlation$=\rho$.
%the grid with mesh $1/h$ on $[0,1]$.
%three, on
%the grid with mesh $3/h$ on $[0,1]$.
%on
%the grid with mesh $6/h$ on $[0,1]$.
%from the grid estimation with mesh $1/h$ on $[0,1]$. We randomly pick
%300 dates for the second Brownian motion from the grid estimation with
%mesh $1/h$ on $[0,2]$. Estimation of the lag parameter
%the grid with mesh $1/h$ on $[0,1]$.
The estimation results are displayed in Table~\ref{tab1} below.
With no surprise, for a given mesh $h_n$, the difficulty of the
estimation problem increases as $\rho$ decreases.

In the fine grid approximation case (FG) with mesh $h_n=10^{-3}$,
the lead-lag parameter $\vartheta$ belongs to ${\mathcal G}^n$ exactly.
Therefore, the contrast ${\mathcal U}^n(\widehat\vartheta)$ is
close to $0$ for all values $\widetilde\vartheta\in{\mathcal
G}^n$, except perhaps for the exact value $\widetilde\vartheta=
\vartheta$. This is illustrated in Figure~\ref{fig1} and Figure~\ref{fig2} below,
where we display the values or ${\mathcal U}^n(\widehat
\vartheta_n)$. Note how more scattered are the values of
${\mathcal U}^n(\widehat\vartheta_n)$ for $\rho=0.25$ compared to
$\rho=0.75$. This is, of course, no surprise.\vadjust{\goodbreak}

%f2 ###
\begin{figure}

\includegraphics{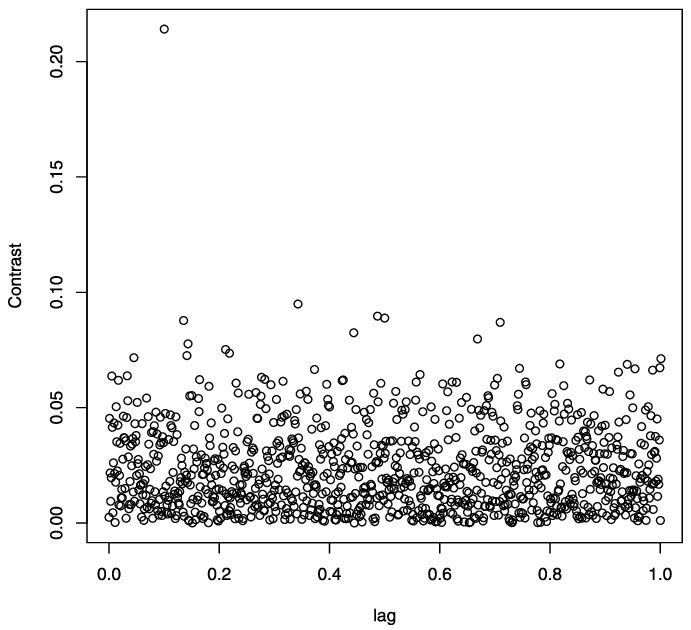}

\caption{Same setting as in Figure \protect\ref{fig1} for $\rho=0.25$.
The value $\max_{\widetilde\vartheta\in{\mathcal G}^n}
| {\mathcal U}^n(\widetilde\vartheta)|$ is still
correctly located.}\label{fig2}
\end{figure}

%f3 ###
\begin{figure}[b]

\includegraphics{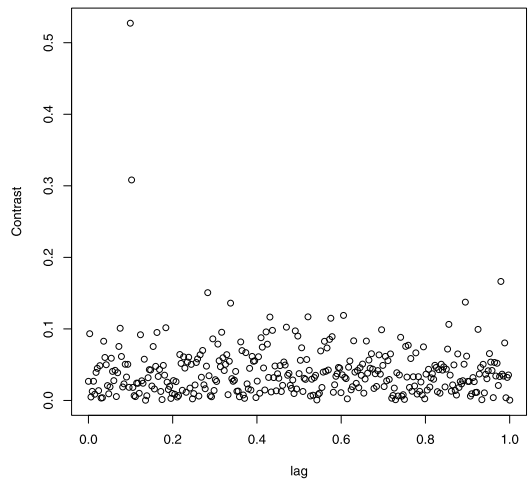}

\caption{Moderate grid case (MG). Over one
simulation: displayed values of $| {\mathcal U}^n(\widetilde
\vartheta)|$ for $\widetilde\vartheta\in{\mathcal G}^n$
with mesh $h_n=10^{-3}$ and $\rho=0.75$. The value
$\max_{\widetilde\vartheta\in{\mathcal G}^n}|{\mathcal
U}^n(\widetilde\vartheta)|$ is still well located. We begin
to see the effect of the maximization over a grid ${\mathcal G}^n$
which does not match exactly with the true value $\vartheta$ with
the appearance of a second maximum.}\label{fig3}
\end{figure}

%f4 ###
\begin{figure}

\includegraphics{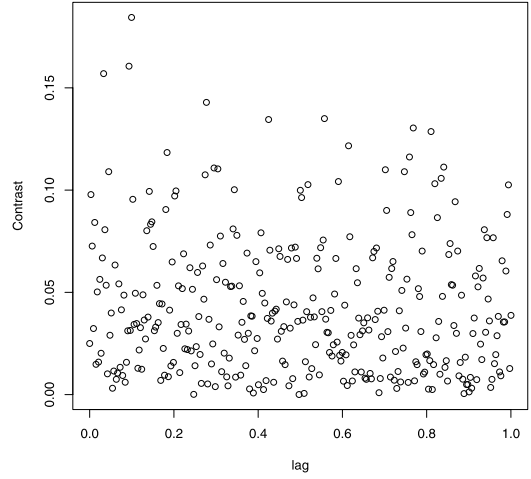}

\caption{Same setting as in Figure \protect\ref{fig3} for $\rho=0.25$.
The value $\max_{\widetilde\vartheta\in{\mathcal G}^n}
|{\mathcal U}^n(\widetilde\vartheta)|$ is still correctly
located, but the overall shape of $|{\mathcal U}^n(\widetilde
\vartheta)|$ deteriorates.}\label{fig4}\vspace*{-4pt}
\end{figure}

%f5 ###
\begin{figure}[b]

\includegraphics{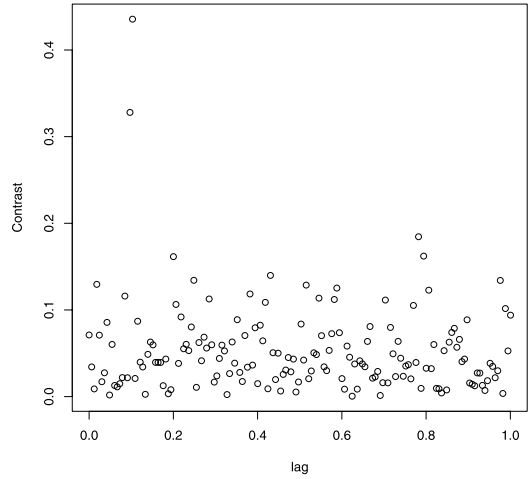}

\caption{Coarse grid case (CG). Over one simulation:
displayed values of $| {\mathcal U}^n(\widetilde
\vartheta)|$ for $\widetilde\vartheta\in{\mathcal G}^n$
with mesh $h_n=10^{-3}$ and $\rho=0.75$. The value
$\max_{\widetilde\vartheta\in{\mathcal G}^n}|{\mathcal
U}^n(\widetilde\vartheta)|$ is still well located. The fact
that ${\mathcal G}^n$ does not match~$\vartheta$ appears more
clearly than in Figure \protect\ref{fig3}.}\label{fig5}
\end{figure}

For the moderate grid (MG) and the coarse grid (CH) cases, the lead-lag
parameter $\vartheta\notin{\mathcal G}^n.$
%does not belong to the grid with mesh
Hence, ${\mathcal U}^n(\widetilde\vartheta)$ is close to $0$ for
almost all values of ${\mathcal G}^n$ except but two.
%The estimation of the correlation parameter
%is "split" between these two "admissible" values.
When $\rho$ is small,
%if the number of data becomes too
%small,
the statistical error in the estimation of $\rho$ is such that
$|\max_{\widetilde\vartheta\in{\mathcal G}^n} {\mathcal
U}^n(\widetilde\vartheta)|$ is not well located\ anymore. The
error in the estimation can then be substantial, but is
nevertheless consistent with our convergence result. This is
illustrated in Figures~\ref{fig3} to~\ref{fig8} below.

When $\rho$ decreases or when the mesh $h_n$ of the grid
increases, the performance of $\widehat\vartheta_n$ deteriorates,
as shown in Figures~\ref{fig7} and~\ref{fig8} below.

%f6 ###
\begin{figure}

\includegraphics{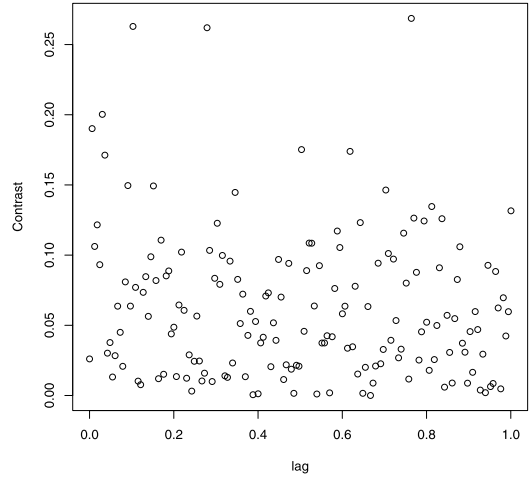}

\caption{Same setting as in Figure \protect\ref{fig5} for $\rho=0.25$.
The value $\max_{\widetilde\vartheta\in{\mathcal G}^n}
|{\mathcal U}^n(\widetilde\vartheta)|$ is no longer
correctly located.}\label{fig6}
\end{figure}

%f7 ###
\begin{figure}[b]

\includegraphics{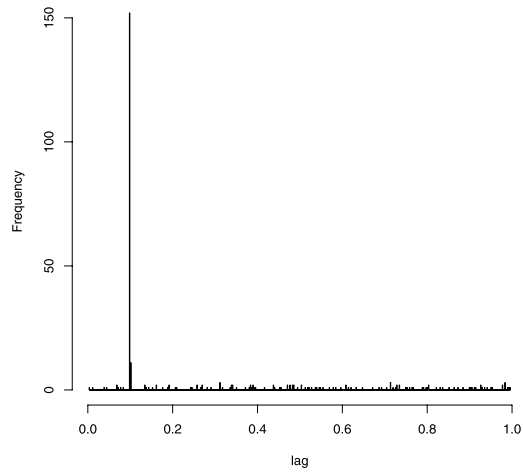}

\caption{Moderate grid case (MG). Histogram of the
values of $\widehat\vartheta_n$ with true value $\vartheta=0.1$
over $300$ simulations for $\rho=0.25$.}\label{fig7}
\end{figure}

%s5.3 ###
\subsection{Non-synchronous data}
We randomly pick 300 sampling times for $X$ over $[0,1]$ uniformly
over a grid of mesh size $10^{-3}$. We randomly pick $300$
sampling times for $Y$ likewise, and independently of the sampling
for $X$. The data generating process is the same as in Section
\ref{DGP}. In Table~\ref{tab2}, we display the estimation results for 300
simulations, in the fine gird case (FG) with $\vartheta=0.1$ and
$\rho=0.75$.

The histograms for the case $\rho=0.5$ and $\rho=0.25$ are
displayed in Figures~\ref{fig9} and~\ref{fig10}.

%f8 ###
\begin{figure}

\includegraphics{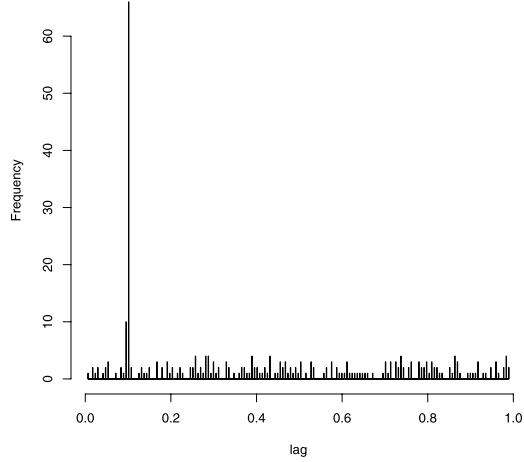}

\caption{Coarse grid case (CG). Histogram of the
values of $\widehat\vartheta_n$ with true value $\vartheta=0.1$
over $300$ simulations for $\rho=0.25$.}\label{fig8}
\end{figure}

%t2 ###
\begin{table}[b]
\caption{Estimation of $\vartheta=0.1$ on 300 simulated
samples for $\rho=0.75$ and non-synchronous data}\label{tab2}
\begin{tabular*}{\textwidth}{@{\extracolsep{\fill}}llllllll@{}}
\hline
% after \  \hline or \cline{col1-col2} \cline{col3-col4} ...
$\widehat\vartheta$&0.099 & 0.1 & 0.101& 0.102 &0.103&0.104&0.105\\
\hline
FG, $\rho=0.75$ &16& 106 & 107 & 46&19&4&2\\
\hline
\end{tabular*}
\end{table}

%s6 ###
\section{A numerical illustration on real data} \label{illustration
real data}

%f9 ###
\begin{figure}

\includegraphics{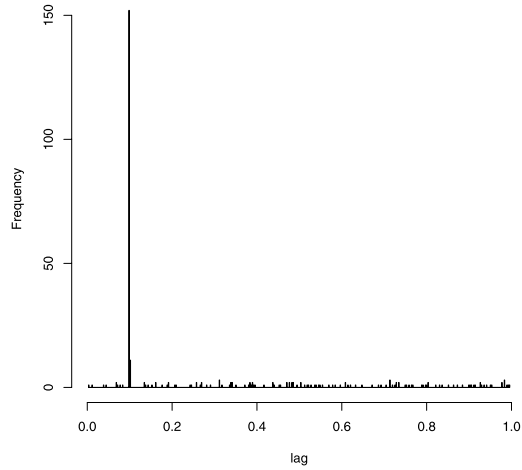}

\caption{Fine grid case (FG), non-synchronous data.
Histogram of the values of $\widehat\vartheta_n$ with true value
$\vartheta=0.1$ over $300$ simulations for
$\rho=0.5$.}\label{fig9}\vspace*{-3pt}
\end{figure}

%f10 ###
\begin{figure}

\includegraphics{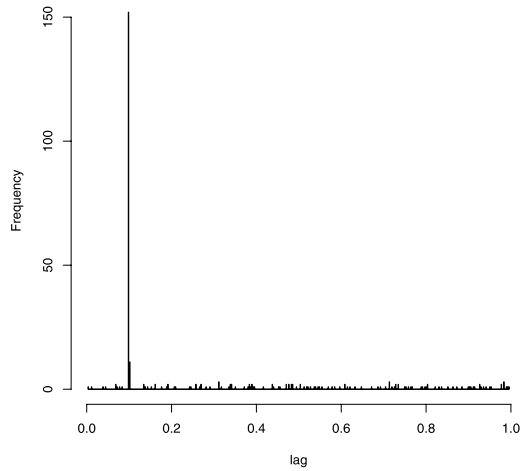}

\caption{Fine grid case (FG), non-synchronous data.
Histogram of the values of $\widehat\vartheta_n$ with true value
$\vartheta=0.1$ over $300$ simulations for $\rho=0.25$. The
performances of $\widehat\vartheta_n$ clearly deteriorates as
compared to Figure \protect\ref{fig9}.}\label{fig10}
\end{figure}

%In this section, we apply our proce

%s6.1 ###
\subsection{The data set}

We study here the lead-lag relationship between the following two
financial assets:

\begin{itemize}[-]
 \item[-] The future contract on the DAX index (FDAX for short),
with maturity December 2010.

 \item[-] The Euro-Bund future contract (Bund for short), with
maturity December 2010, which is an interest rate product based on
a notional long-term debt instruments issued by the Federal
Republic of Germany.
\end{itemize}

These two assets are electronically traded on the EUREX market,
and are known to be highly liquid. Our data set has been provided
by the company QuantHouse
EUROPE/ASIA\footnote{\url{http://www.quanthouse.com}.}. It consists in
all the trades for 20 days of October 2010. Each trading day
starts at 8.00 am CET and finishes at 22.00 CET, and the accuracy
in the timestamp values is one millisecond.\vadjust{\goodbreak}

%
%f11 ###
\begin{figure}

\includegraphics{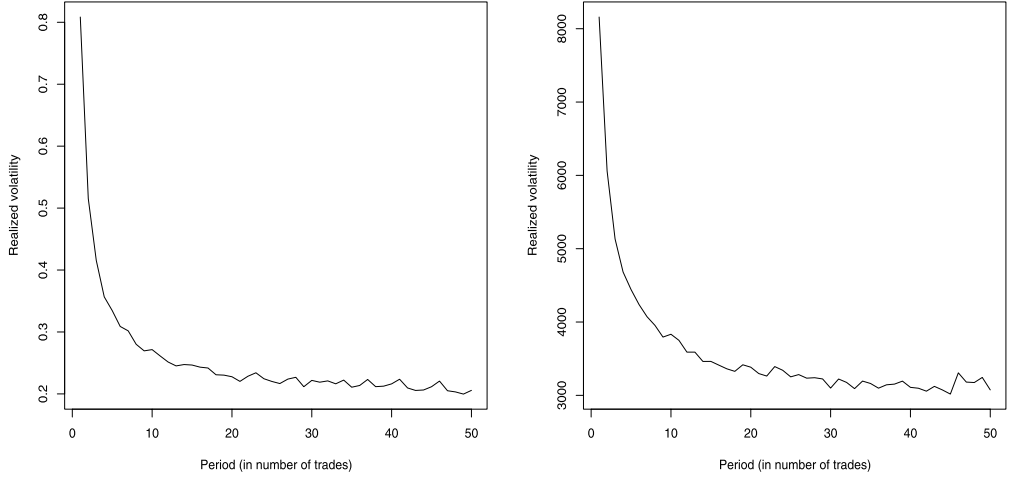}

\caption{Signature plot for the Bund (left) and the
FDAX (right) for 2010, October 13.}\label{fig11}
\end{figure}

%s6.2 ###
\subsection{Methodology: A one day analysis}

In order to explain our methodology, we take the example of a
representative day: 2010, October~13.

\subsubsection*{Microstructure noise}
Since high-frequency data are concerned, we need to incorporate
microstructure noise effects, at least at an empirical level. A
classical way to study the intensity of the microstructure noise is
to draw the signature plot (here in trading time). The signature
plot is a function from $\mathbb{N}$ to $\mathbb{R}^+$. To a given
integer $k$, it associates the sum of the squared increments of the
traded price (the realized volatility) when only 1 trade out of $k$
is considered for computing the traded price. If the price were
coming from a continuous-time semi-martingale, the signature plot
should be approximately flat. In practice, it is decreasing, as
shown by Figure~\ref{fig11}.

%f12 ###
\begin{figure}

\includegraphics{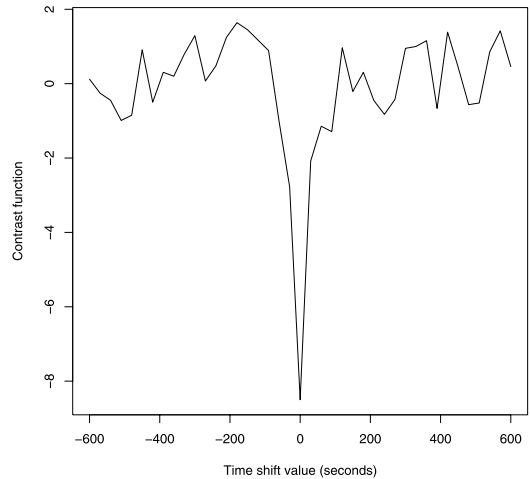}

\caption{The function ${\mathcal U}^n$ for 2010, October
13, time
shift values between $-10$ minutes and $10$ minutes, on a grid
with mesh 30 seconds. The contrast is obtained by taking the absolute
value of ${\mathcal U}^n$.}\label{fig12}
\end{figure}

According to Figure~\ref{fig11}, for all our considered day, we subsample our
data so that we keep one trade out of 20. On 2010, October 13, after
subsampling, it remains 2018 trades for the Bund and 3037
trades for the FDAX.

\subsubsection*{Construction of the contrast function}

The second step is to compute our contrast function. Here the Bund
plays the role of $X$ and the FDAX the role of $Y$. Therefore, if
the estimated value is positive, it means that the Bund is the
leader asset and the FDAX the lagger asset, and conversely. To
have a first idea of the lead-lag value, we consider our contrast
function for a time shift between $-10$ minutes and $10$ minutes,
on a grid with mesh 30 seconds. The result of this computation for
October 2010, 13 is given in Figure~\ref{fig12}.

\begin{figure}

\includegraphics{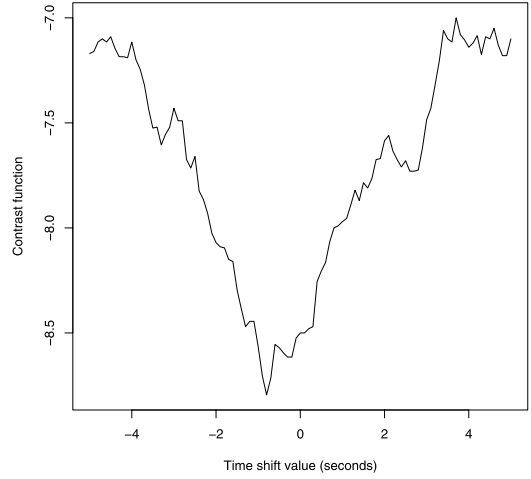}

\caption{The function ${\mathcal U}^n$ for 2010, October
13, time
shift values between $-5$ seconds and $5$ seconds, on a grid with
mesh 0.1 second. The contrast is obtained by taking the absolute value
of ${\mathcal U}^n$.}\label{fig13}
\end{figure}

From Figure~\ref{fig12}, we see that the lead-lag value is close to zero.
Thus, we then compute the contrast function for a time shift
between $-5$ seconds and $5$ seconds, on a grid with mesh 0.1
second. The result of this computation for 2010, October 13 is
given in Figure~\ref{fig13}.

%f13 ###

From Figure~\ref{fig13}, we can conclude that on 2010, October 13, the FDAX
seems to lead the Bund, with a small lead lag value of $-0.8$
second.
%s6.3 ###
\subsection{Systematic results over a one-month period}
We now give, in Figure~\ref{fig14}, the results for all the days of October
2010.
%

%f14 ###
\begin{figure}
\fontsize{9}{11}{\selectfont{
\begin{tabular*}{\textwidth}{@{\extracolsep{\fill}}llll@{}}
\hline
 & Number of trades for the& Number of trades for the& Lead-lag \\
Day& bund (after subsampling) & FDAX (after subsampling) &
(seconds)\\
\hline
1 October 2010 & 2847 & 4215 & $-0.2$ \\
5 October 2010 & 2213 & 3302 & $-1.1$ \\
6 October 2010 & 2244 & 2678 & $-0.1$ \\[3pt]
7 October 2010 & 1897 & 3121 & $-0.5$ \\
8 October 2010 & 2545 & 2852 & $-0.6$ \\
11 October 2010 & 1050 & 1497 & $-1.4$ \\[3pt]
12 October 2010 & 2265 & 3018 & $-0.8$ \\
13 October 2010 & 2018 & 3037 & $-0.8$ \\
14 October 2010 & 2057 & 2625 & $-0.0$ \\[3pt]
15 October 2010 & 2571 & 3269 & $-0.7$ \\
18 October 2010 & 1727 & 2326 & $-2.1$ \\
19 October 2010 & 2527 & 3162 & $-1.6$ \\[3pt]
20 October 2010 & 2328 & 2554 & $-0.5$ \\
21 October 2010 & 2263 & 3128 &$ -0.1$ \\
22 October 2010 & 1894 & 1784 & $-1.2$ \\[3pt]
25 October 2010 & 1501 & 2065 & $-0.4$ \\
26 October 2010 & 2049 & 2462 & $-0.1$ \\
27 October 2010 & 2606 & 2864 & $-0.6$ \\[3pt]
28 October 2010 & 1980 & 2632 & $-1.3$ \\
29 October 2010 & 2262 & 2346 & $-1.6$ \\
\hline
\end{tabular*}}}
\caption{Estimated lead-lag values for October
2010.}\label{fig14}
\end{figure}

The results of Figure~\ref{fig14} seem to indicate that, on average, the FDAX
tends to
lead the Bund. Indeed, the estimated lead-lag values are
systematically negative. Of course these results have to be taken
with care since the estimated values are relatively small (the order of
one second); however, dealing with highly traded assets
on electronic markets, the order of magnitude of the lead-lag values
that we find are no surprise and are consistent with common knowledge.
A~possible interpretation -- yet speculative at the exploratory level
intended here -- for the presence of such lead-lag effects is the
difference between the tick sizes of the different assets. Indeed, the
negative values could mean that the tick size of the FDAX can be
considered smaller than those of the Bund.

\begin{appendix}

\section*{Appendix}\label{appendix}

%(a): %It is straightforward from definition.
%For $t\in[0,\vartheta)$, $\{\underline{I^n_\vartheta}\leq t\}=\varnothing
%$\{\underline{I_\vartheta}\leq t\}\in\in\calf^{\vartheta^*+
%by definition. The same reasoning is valid for
%$\overline{I^n_\vartheta}$.
%&=&
%=
%
\setcounter{subsection}{0}
%s6.1 ###
\subsection{\texorpdfstring{Proof of Proposition \protect\ref{asympt_behaviour_contrast}}{Proof of Proposition 1}}
For notational clarity, for a
given interval $I = (\underline{I}, \overline{I}]$, we may sometimes
write $X(\underline{I},\overline{I})$ instead of $X(I)$ when no
confusion is possible. In the Bachelier case with lead-lag parameter
$\vartheta\in\Theta$, we work with the following explicit
representation of the observation process:
\renewcommand{\theequation}{\arabic{equation}}
\setcounter{equation}{12}
\begin{equation} %\label{leadlagrepresentation}
\cases{
X_t  =  x_0+\sigma_1 B_t,\vspace*{2pt}\cr
Y_t =  y_0+ \sigma_2 \bigl(\rho B_{t-\vartheta}+\sqrt{1-\rho^2}
W_{t-\vartheta}\bigr),}
\end{equation}
where $B$ and $W$ are two independent Brownian motions. We have
\[
{\mathcal U}^n(\tilde{\vartheta}) = \sum_{0 \leq i\Delta_n \leq T}
X
\bigl((i-1)\Delta_n, i\Delta_n\bigr)\tau_{-\tilde{\vartheta}}Y
\bigl((i-1)\Delta
_n, i\Delta_n\bigr)
= \sigma_1\sigma_2 \sum_{0 \leq i\Delta_n \leq T} \chi_i^n(\tilde
{\vartheta}),
\]
with
\[
\chi_i^n(\tilde{\vartheta})=B\bigl((i-1)\Delta_n, i\Delta_n\bigr)
\bigl[\rho  \tau_{\vartheta-\tilde{\vartheta}}B\bigl((i-1)\Delta_n,
i\Delta
_n\bigr)+\sqrt{1-\rho^2} \tau_{\vartheta-\tilde{\vartheta}}W
\bigl((i-1)\Delta_n, i\Delta_n\bigr)\bigr].
\]
We have
\[
\E[\chi_i^n(\tilde{\vartheta})] = \rho  \E\bigl[B
\bigl((i-1)\Delta_n, i\Delta_n\bigr)\tau_{\vartheta-\tilde{\vartheta
}}B
\bigl((i-1)\Delta_n, i\Delta_n\bigr)\bigr]
= \rho \Delta_n  \varphi\bigl(\Delta_n^{-1}(\tilde{\vartheta
}-\vartheta)\bigr),
\]
where $\varphi(x)=(1-|x|)1_{|x|\leq1}$ is the usual hat function.
Assuming further, with no loss of generality, that $T/\Delta_n$ is an
integer, we obtain the representation
\[
{\mathcal U}^n(\tilde{\vartheta}) = \sigma_1\sigma_2T \biggl( \rho
\varphi\bigl(\Delta_n^{-1}(\vartheta-\tilde{\vartheta})\bigr) +
T^{-1}\sum
_{0 \leq i\Delta_n \leq T} \bigl(\chi_i^n(\tilde{\vartheta})-\E
[\chi
_i^n(\tilde{\vartheta})]\bigr) \biggr).
\]
We now assume without loss of generality that $0 \leq\vartheta-\tilde
{\vartheta} \leq\Delta_n$ (the symmetric case being treated the same
way). The sequence of random variables $\chi_i^n(\tilde{\vartheta})$ is
stationary. Moreover, since the random variable $\chi_i^n(\tilde
{\vartheta})$ involves increments of
$W$ and $B$ over a domain included in $[(i-2)\Delta_n,i\Delta_n]$
because $|\vartheta-\tilde{\vartheta}| \leq\Delta_n$, it follows that
$\chi_i^n(\tilde{\vartheta})$ and $\chi_j^n(\tilde{\vartheta})$ are
independent as soon as $|i-j|\geq2$. Moreover, we claim that
\begin{equation} \label{cov nulle}
\operatorname{Cov}(\chi_i^n(\tilde{\vartheta}),\chi_j^n(\tilde
{\vartheta
}))=0  \qquad\mbox{if }  |i-j|=1.
\end{equation}
Therefore, by the central limit theorem, we have that
\[
\Delta_n^{1/2}T^{-1/2}\sum_{0 \leq i\Delta_n \leq T} \bigl(\chi
_i^n(\tilde{\vartheta})-\E[\chi_i^n(\tilde{\vartheta})]\bigr)
\]
is approximately centred Gaussian, with variance
\[
\operatorname{Var}(\chi_1^n(\tilde{\vartheta})).
\]
\subsubsection*{Computation of $\operatorname{Var}(\chi_1^n(\tilde
{\vartheta}))$} To that end, we need to evaluate
\[
\mathrm{I} = \rho^2 \E[(B(0, \Delta_n)\tau_{\vartheta
-\tilde
{\vartheta}}B(0, \Delta_n))^2],
\]
and
\[
\mathrm{II} = (1-\rho^2)\E[(B(0,\Delta_n)\tau_{\vartheta-\tilde
{\vartheta
}}W(0,\Delta_n))^2],
\]
since $B(0, \Delta_n)\tau_{\vartheta-\tilde{\vartheta
}}B(0,
\Delta_n)$ and $B(0, \Delta_n)\tau_{\vartheta-\tilde
{\vartheta}}W(0, \Delta_n)$ are uncorrelated. Writing
\begin{eqnarray*}
&& B(0, \Delta_n)\tau_{\vartheta-\tilde{\vartheta}}B(0,
\Delta_n) \\
&&\quad=  \bigl(B(0,\vartheta-\tilde{\vartheta})+B(\vartheta-\tilde
{\vartheta
},\Delta_n)\bigr)\bigl(B(\vartheta-\tilde{\vartheta}, \Delta
_n)+B(\Delta
_n,\vartheta-\tilde{\vartheta}+\Delta_n)\bigr),
\end{eqnarray*}
taking square and expectation, we readily obtain that
\begin{eqnarray*}
\mathrm{I} & =& 2 \rho^2(\vartheta-\tilde{\vartheta})
\bigl(\Delta
_n- (\vartheta-\tilde{\vartheta})\bigr)+\rho^2(\vartheta-\tilde
{\vartheta
})^2 + 3\rho^2\bigl(\Delta_n- (\vartheta-\tilde{\vartheta})
\bigr)^2
\\
& =&\rho^2\bigl(\Delta_n^{2}\bigl(1+2\varphi\bigl(\Delta
_n^{-1}(\vartheta
-\tilde{\vartheta})\bigr)^2\bigr)\bigr).
%)\varphi(\Delta_n^{-1}(\vartheta_0-\vartheta)) + {\mathcal
%O}(\Delta_n^2) \\
%& = 2\Delta_n \varphi(\Delta_n^{-1}(\vartheta_0-\vartheta)) + {
\end{eqnarray*}
%
%By symmetry, the same computation is valid for $-\Delta_n \leq
Concerning II, since $B$ and $W$ are independent, we readily have
\[
\mathrm{II} = (1-\rho^2)\Delta_n^2,
\]
therefore, from $\E[\chi_1^n(\tilde{\vartheta})]= \rho \Delta_n
\varphi(\Delta_n^{-1}(\tilde{\vartheta}-\vartheta))$, we
finally infer
\[
\Delta_n^{-2}\operatorname{Var}(\chi_1^n(\tilde{\vartheta})) =
1+\rho
^2\varphi\bigl(\Delta_n^{-1}(\vartheta-\tilde{\vartheta})\bigr)^2
% \end{equation}
\]
from which Proposition~\ref{asympt_behaviour_contrast} follows. It
remains to prove \eqref{cov nulle}. By stationarity, this amounts to evaluate
\[
\rho^2\E[B(0,\Delta_n)B(\vartheta-\tilde{\vartheta}, \Delta
_n+\vartheta-\tilde{\vartheta})B(\Delta_n,2\Delta_n)B(\Delta
_n+\vartheta
-\tilde{\vartheta}, 2\Delta_n+\vartheta-\tilde{\vartheta})
]-\E[\chi
_1^n(\vartheta)]^2.
\]
To that end, we split each of the terms as follows:
\begin{eqnarray*}
B(0,\Delta_n) & =& B(0,\vartheta-\tilde{\vartheta})+B(\vartheta
-\tilde
{\vartheta},\Delta_n),\\
B(\vartheta-\tilde{\vartheta}, \Delta_n+\vartheta-\tilde
{\vartheta}) &
= &B(\vartheta-\tilde{\vartheta}, \Delta_n)+B(\Delta_n+\vartheta
-\tilde
{\vartheta}),\\
B(\Delta_n,2\Delta_n) & =& B(\Delta_n,\Delta_n+\vartheta-\tilde
{\vartheta
}) + B(\Delta_n+\vartheta-\tilde{\vartheta},2\Delta_n),\\
B(\Delta_n+\vartheta-\tilde{\vartheta}, 2\Delta_n+\vartheta
-\tilde
{\vartheta})& = &B(\Delta_n+\vartheta-\tilde{\vartheta}, 2\Delta
_n) +
B(2\Delta_n, 2\Delta_n+\vartheta-\tilde{\vartheta}).
\end{eqnarray*}
Using the stochastic independence of each of these terms, multiplying
and integrating, we easily obtain
\begin{eqnarray*}
&& \rho^2\E[B(0,\Delta_n)B(\vartheta-\tilde{\vartheta},
\Delta
_n+\vartheta-\tilde{\vartheta})B(\Delta_n,2\Delta_n)B(\Delta
_n+\vartheta
-\tilde{\vartheta}, 2\Delta_n+\vartheta-\tilde{\vartheta})]
\\
&&\quad=  \rho^2\Delta_n^2 \varphi\bigl(\Delta_n^{-1}(\vartheta-\tilde
{\vartheta})\bigr)^2 = \E[\chi_1^n(\tilde{\vartheta})]^2.
\end{eqnarray*}

%s6.2 ###
\subsection{\texorpdfstring{Proof of Proposition \protect\ref{noclt}}{Proof of Proposition 2}}
Suppose that $\Delta_n^{-1}(\widehat\vartheta_n-\vartheta
)\rightarrow
Z$, in law, for some random random variable $Z$. For $a\in\mathbb{R}$,
we write $a^{[ n]},$ the best approximation of $a$ by a point of the
form $k\Delta_n$, $k\in\mathbb{Z}$ and $a^{\lfloor n\rfloor}$, the best
approximation of $a$ by a point smaller or equal to $a$ and of the form
$k\Delta_n$, $k\in\mathbb{Z}$. We have
\[
\Delta_n^{-1}(\widehat\vartheta_n-\vartheta)=\Delta
_n^{-1}\bigl(\widehat
\vartheta_n-\widehat\vartheta_n^{\lfloor n\rfloor}\bigr)+\Delta
_n^{-1}\bigl(\widehat\vartheta_n^{\lfloor n\rfloor}-\vartheta\bigr).
\]
The first term in the right-hand side of the equality is smaller than
$\Delta_n^{-1}h_n$ and so converges to zero. The second term can be
written as
\[
\Delta_n^{-1}\bigl(\widehat\vartheta_n^{\lfloor n\rfloor}-\vartheta
^{[n]}\bigr)+\Delta_n^{-1}\bigl(\vartheta^{[n]}-\vartheta\bigr)=T_{1,n}+T_{2,n},
\]
say. The sequence $T_{1,n}$ is a random sequence of integers, and
$T_{2,n}$ is a deterministic sequence with values in $[0,1/2]$ which
does not converge. Let $\psi_n$ be a subsequence such that $T_{2,\psi
_n}\rightarrow l$ with $l\in(0,1/2]$. Then $T_{1,\psi_n}$ converges in
law to $Z-l$ which implies that the support of $Z$ is included in $\{
z+l,z\in\mathbb{Z}\}$. Consider now $\widetilde\psi_n$ such that
$T_{2,\widetilde\psi_n}\rightarrow l'$ with $l'\in[0,1/2]$, $l'\neq
l$. In the same way, we get that the support of $Z$ is also included in
$\{z+l',z\in\mathbb{Z}\}$, a contradiction.

%s6.3 ###
\subsection{\texorpdfstring{Proof of Lemma \protect\ref{lagest-1}}{Proof of Lemma 1}}

\subsubsection*{Preliminary results} We first prove the following
results.
\begin{lemm}\label{lagest-1bis} Work under Assumption \textup{B2}, under the
slightly more general assumption that for all $I =[\underline{I,
\overline{I}}) \in{\mathcal I}$, the random variables $\underline{I}$
and $\overline{I}$ are $\F$-stopping times.
%Assume Condition $[A1^\flat]$. %, $[A2]$ and $[A3]$.
%
\begin{enumerate}[(a)]
\item[(a)]
%If $\vartheta\geq\vartheta^*+\bar{r}_n$, then
%$\underline{I^n_\vartheta}$ and $\overline{I^n_\vartheta}$ are
%$\F^{\vartheta^*+\bar{r}_n}$-stopping times.
%
If $\widetilde\vartheta\geq\vartheta+v_n$, then
for any $\F$-stopping time $\sigma$ and $t\in\bbR_+$, $\sigma
+\widetilde\vartheta$ is an $\F^{\vartheta+v_n}$-stopping time. In
particular, the random variables $\underline{I^n_{\widetilde\vartheta}}$
and $\overline{I^n_{\widetilde\vartheta}}$ are $\F^{\vartheta
+v_n}$-stopping times.

\item[(b)]
For each $J \in{\mathcal J}$, we have
%, where $J^n_t=J^n\cap[0,t]$.
%As a consequence of (ii), we will have
\[\label{d081021-1}
\calf^{\vartheta+v_n}_{\overline{J^n}}
\subset
\calf^{\vartheta+v_n}_{\underline{J^n}+v_n}
=
\calf^{\vartheta}_{\underline{J^n}},
\]
and for each $I \in{\mathcal I}$,
\[\label{d081021-2}
\calf_{\overline{I^n}}
=
\calf^{\vartheta+v_n}
_{\overline{I^n_{(\vartheta+v_n)}}}.
\]

\item[(c)]
Suppose that $\widetilde\vartheta\geq\vartheta+\ep_n$
and $2v_n\leq\ep_n$. Then for any random variable $X'$ measurable
w.r.t. $\calf_{\overline{I^n}}$, the random variables
$X'1_{\{\underline{I^n_{\widetilde\vartheta}}\leq\overline{J^n}\}}$
and $X'1_{\{\underline{I^n_{\widetilde\vartheta}}<\overline{J^n}\}}$
are $\calf^{\vartheta}_{\underline{J^n}}$-measurable.
\end{enumerate}
\end{lemm}

\begin{pf}
\textit{Proof of} (a). For any $\F$-stopping time $\sigma$ and
$t\in\bbR_+$,
\begin{eqnarray*}
\{\sigma+\widetilde\vartheta\leq t \}
&=&\{\sigma\leq t-\widetilde\vartheta\}
=
\bigl\{\sigma\leq\bigl(t-(\widetilde\vartheta-\vartheta-v_n)\bigr)-\vartheta-v_n
\bigr\}
\\
&\in&
\calf^{\vartheta+v_n}_{t-(\widetilde\vartheta-\vartheta-v_n)}
\subset\calf^{\vartheta+v_n}_t.
\end{eqnarray*}

\textit{Proof of} (b). Note first that under Assumption B2,
the $\F^{\vartheta+v_n}$-stopping time $\underline{J^n}$ is
in particular an $\F^{\vartheta}$-stopping time;
thus $\calf^{\vartheta}_{\underline{J^n}}$ is a $\sigma$-field.
%Furthermore we assume $2\bar{r}_n\leq\ep_n$ for all $n\in\bbN$.
Moreover, since $\overline{J^n}$ and $\underline{J^n}+v_n$
are $\F^{\vartheta+v_n}$-stopping times by definition,
both
$\calf^{\vartheta+v_n}_{\overline{J^n}}$ and
$\calf^{\vartheta+v_n}_{\underline{J^n}+v_n}$ are
$\sigma$-fields, and also the inclusion is trivial from
$\overline{J^n}\leq\underline{J^n}+v_n$.
To obtain the equality, it suffices to observe
that each of the conditions
``${\mathcal A}\in\calf^{\vartheta+v_n}_{\underline{J^n}+v_n}$''
and ``${\mathcal A}\in\calf^{\vartheta}_{\underline{J^n}}$''
is equivalent to the condition
\[
{\mathcal A}\cap\{\underline{J^n}\leq t-v_n\}
\in\calf^{\vartheta}_{t-v_n}
\]
for all $t\in\bbR_+$. The second equality is proved in the same way.
%%Recall that $\vartheta^*\geq0$.

\textit{Proof of} (c). Since $\overline{J^n}$ and $\underline
{I^n_{\widetilde\vartheta}}$ are $\F^{\vartheta+v_n}$-stopping times
by assumption, we have
\[
\{\underline{I^n_{\widetilde\vartheta}}\leq\overline{J^n}\}\in
\calf
^{\vartheta+v_n}
_{\overline{J^n}}\subset\calf^{\vartheta}_{\underline{J^n}},
\]
the last inclusion following from (b).
%{From} definition of $\underlin{I^n}$ and $\overline{I^n}$,
If $\underline{I^n_{\widetilde\vartheta}}\leq\overline{J^n}$, then
\[
\overline{I^n}
\leq\underline{I^n}+v_n
\leq\overline{J^n}-\widetilde\vartheta+v_n
\leq\overline{J^n}-\vartheta-v_n,
\]
which implies $\overline{I^n_{ \vartheta+v_n}}\leq\overline{J^n}.$
Thus
\[
X'1_{\{\underline{I^n_{\widetilde\vartheta}}\leq\overline{J^n}\}}
=
X'
1_{\{\overline{I^n_{(\vartheta+v_n)}}\leq\overline{J^n}\}}\times
1_{\{\underline{I^n_{\widetilde\vartheta}}\leq\overline{J^n}\}}.
\]
We have that $X'$ is measurable with respect to
$\calf_{\overline{I^n}}=\calf^{\vartheta+v_n}
_{\overline{I^n_{(\vartheta+v_n)}}}$. Also
%_{\overline{I^n}+\vartheta^*+\bar{r}_n}$, and also
$\overline{I^n_{(\vartheta+v_n)}}$ is
a stopping time with respect to
$\F^{\vartheta+v_n}$ by (a). Consequently,
$X'1_{\{\overline{I^n_{(\vartheta+v_n)}}\leq\overline{J^n}\}}$
is $\calf^{\vartheta+v_n}_{\overline{J^n}}$-measurable, hence
$\calf^{\vartheta}_{\underline{J^n}}$-measurable.
Eventually, $X'1_{\{\underline{I^n_{\widetilde\vartheta}}\leq
\overline
{J^n}\}}$
is $\calf^{\vartheta}_{\underline{J^n}}$-measurable. The other
statement is proved the same way.
\end{pf}

\begin{pf*}{Proof of Lemma \protect\ref{lagest-1}}
We have
\[
X'K(I^n_{\widetilde\vartheta},J^n)
=
X'1_{\{\underline{I^n_{\widetilde\vartheta}}\leq\underline{J^n}\}}
1_{\{
\underline{J^n}
<\overline{I^n_{\widetilde\vartheta}}\}}
+
X'1_{\{\underline{I^n_{\widetilde\vartheta}}>\underline{J^n}\}}
1_{\{\overline{J^n}> \underline{I^n_{\widetilde\vartheta}}\}}.
\]
Since $\widetilde\vartheta\geq\vartheta+\ep_n \geq\vartheta+v_n$,
%both $\underline{I^n_\vartheta}$ and $\underline{J^n}$ are
both $\underline{I^n_{\widetilde\vartheta}}$ and $\overline
{I^n_{\widetilde\vartheta}}$ are
$\F^{\vartheta+v_n}$-stopping times. Therefore,
the second term on the right-hand side of the above equality
is $\calf^{\vartheta}_{\underline{J^n}}$-measurable
by (c) of Lemma~\ref{lagest-1bis}.
%Measurability of the first term is easier.

Now we notice that
if
$\underline{I^n_{\widetilde\vartheta}}\leq\underline{J^n}$, then
$\underline{I^n_{\widetilde\vartheta}}\leq\overline{J^n}$, therefore
\[
X'1_{\{\underline{I^n_{\widetilde\vartheta}}\leq\underline{J^n}\}}
1_{\{\underline{J^n}<\overline{I^n_{\widetilde\vartheta}}\}}
=
\bigl(X'1_{\{\underline{I^n_{\widetilde\vartheta}}\leq\overline
{J^n}\}
}\bigr)
\times\bigl(1_{\{\underline{I^n_{\widetilde\vartheta}}\leq
\underline
{J^n}\}}
1_{\{\underline{J^n}<\overline{I^n_{\widetilde\vartheta}}\}}\bigr).
\]
The first factor on the right-hand side of the above
equality is $\calf^{\vartheta}_{\underline{J^n}}$-measurable
by (c) of Lemma~\ref{lagest-1bis}, and the second factor is
obviously
$\calf^{\vartheta}_{\underline{J^n}}$-measurable. This completes the proof.
%After all, we obtain the desired measurability.
\end{pf*}

%s6.4 ###
\subsection{\texorpdfstring{Proof of Lemma \protect\ref{stoppingtime2}}{Proof of Lemma 2}}

Let us fix $I \in{\mathcal I}$. Let
\[
T_J = \cases{
\overline{I^n}-v_n &\quad$\mbox{on }  \{\overline{\circj_{-\delta_n}}
>\overline{I^n}\} $,\vspace*{2pt}\cr
\overline{\circj}_{-\delta_n} &\quad $\mbox{on }  \{\overline
{\circj
_{-\delta_n}} \leq\overline{I^n}\}.$}
\]
%
%Let $\G_n=(\calf_{u-\bar{r}_n})_{u\in[-\bar{\vartheta}+\bar{r}_n,
We know that $\overline{I^n}-v_n$ is an $\F$-stopping time by
Assumption B2, and also that $\overline{\circj_{-\delta_n}}-v_n$
is an $\F$-stopping time due to $\delta_n\leq0$. Let us show
first that the $T_J$s are $\F$-stopping times. Let $t \in
[-\delta, T+\delta]$. Let
\[
{\mathcal A}_1 = \{
\overline{I^n}-v_n \leq t,  \overline{\circj_{-\delta_n}}-v_n
>\overline{I^n}-v_n
\}
\]
and
\[
{\mathcal A}_2 = \{
\overline{\circj}_{-\delta_n} \leq t,
\overline{\circj_{-\delta_n}}-v_n \leq\overline{I^n}-v_n \}.
\]
It is obvious that ${\mathcal A}_1 \in\calf_t$ since
$\overline{I^n}-v_n$ is an $\F$-stopping time and also
\[
\{\overline{\circj_{-\delta_n}}-v_n
\geq\overline{I^n}-v_n\} \in\calf_{\overline{I^n}-v_n}.
\]
For
the term ${\mathcal A}_2$, if $t\in[-\delta,-\delta+v_n]$, then
${\mathcal A}_2=\varnothing\in\calf_{-\delta}\subset\calf_t$.
Otherwise, if $t\in(-\delta+v_n,T+\delta]$, then
\[
{\mathcal
A}_2 = \{ \overline{\circj_{-\delta_n}}-v_n\leq t-v_n,
\overline{\circj_{-\delta_n}}-v_n \leq\overline{I^n}-v_n \}
\in\calf_{t-v_n}\subset\calf_t.
\]
Eventually, we have
$\{T_J\leq t\}\in\calf_t$; hence $T_J$ is an $\F$-stopping time.

In conclusion, there exists at least one
$\overline{\circj_{-\delta_n}}$ in
$[\overline{I^n}-v_n,\overline{I^n}]$. Therefore, we have
$M^I=\sup_JT_J$, and this implies that $M^I$ is also an
$\F$-stopping time.
\end{appendix}

\section*{Acknowledgements} This work was originated by
discussions between M. Hoffmann and M. Rosenbaum with S.~Pastukhov
from the Electronic Trading Group research team of P. Gu\'evel at
BNP-Paribas. We are grateful to M. Musiela, head of the Fixed
Income research at BNP-Paribas, for his constant support and
encouragements. We also thank E. Bacry, K. Al Dayri and Tuan Nguyen
for inspiring discussions. Japan Science and Technology supported
the theoretical studies in this work. N.~Yoshida's research was
also supported by Grants-in-Aid for Scientific Research No.~19340021,
the global COE program, ``The research and training
center for new development in mathematics'' of Graduate School of
Mathematical Sciences, University of Tokyo, JST Basic Research
Programs PRESTO and by Cooperative Reserch Program of the
Institue of Statistical Mathematics.

We are grateful to the comments and inputs of two referees and an
associate editor, that help to improve a former version of this work.

% imsref loaded by akundreckaite, 2012-01-04 10:02:03
%

\printhistory

\end{document}